\documentclass[aos,seceqn,nameyear,dvips]{arximspdf}
\usepackage{graphicx}

\doi{10.1214/10-AOS855}
\volume{39}
\issue{2}
\pubyear{2011}
\firstpage{772}
\lastpage{802}

\makeatletter

\makeatletter
\newcommand{\xrightarrow}[1]{\mathrel{\vbox{\m@th\ialign{##\crcr $\hfil\scriptstyle\ #1\ \
\hfil$\crcr\noalign{\kern0.5pt\nointerlineskip}%
\rightarrowfill\crcr}}}}
\makeatother

\newtheorem{satz}{Satz}[section]
\newtheorem{theorem}[satz]{Theorem}
\newtheorem{proposition}[satz]{Proposition}
\newtheorem{corollary}[satz]{Corollary}
\newproclaim{definition}[satz]{Definition}
\newproclaim{remark}[satz]{Remark}

\newcommand{\E}{\mathbb{E}}
\newcommand{\R}{\mathbb{R}}

\newcommand{\N}{\mathbb{N}}
\newcommand{\PP}{\mathbb{P}}
\newcommand{\QQ}{\mathbb{Q}}
\newcommand{\supp}{\operatorname{supp}}
\newcommand{\Id}{\operatorname{Id}}
\newcommand{\Var}{\operatorname{Var}}

\newcommand{\eps}{\varepsilon}
\renewcommand{\phi}{\varphi}
\renewcommand{\theta}{\vartheta}
\newcommand{\abs}[1]{\vert #1 \vert}
\newcommand{\norm}[1]{\Vert #1 \Vert}

\newcommand{\scapro}[2]{\langle #1,#2 \rangle}
\newcommand{\floor}[1]{\lfloor #1 \rfloor}

\newcommand{\eqref}[1]{(\ref{#1})}

\newcommand{\cit}[1]{\citet{#1}}

\makeatother

\begin{document}
\begin{frontmatter}

\title{Asymptotic equivalence for inference on the volatility from noisy observations}
\runtitle{Asymptotic equivalence for inference on the volatility}

\begin{aug}
\author[A]{\fnms{Markus} \snm{Rei{\fontsize{8.8}{8.8}\selectfont\normalfont{\ss}}}\corref{}%
\ead[label=e1]{mreiss@mathematik.hu-berlin.de}}
\runauthor{M. Rei\normalfont{\ss}}
\affiliation{Humboldt-Universit\"{a}t zu Berlin}
\address[A]{Institut f\"{u}r Mathematik\\
Humboldt-Universit\"{a}t zu Berlin\\
Unter den Linden 6\\
D-10099 Berlin\\
Germany\\
\printead{e1}} %adresu isvedimo komanda gale!
\end{aug}

% HISTORY:
\received{\smonth{1} \syear{2010}}
\revised{\smonth{9} \syear{2010}}

% ABSTRACT
%
\begin{abstract}
We consider discrete-time observations of a continuous martingale under
measurement error. This serves as a fundamental model for
high-frequency data in finance, where an efficient price process is
observed under microstructure noise. It is shown that this
nonparametric model is in Le Cam's sense asymptotically equivalent to a
Gaussian shift experiment in terms of the square root of the volatility
function $\sigma$ and a nonstandard noise level. As an application, new
rate-optimal estimators of the volatility function and simple efficient
estimators of the integrated volatility are constructed.
\end{abstract}

% KEYWORDS
%

\end{frontmatter}

%s1 ###
\section{Introduction}

In recent years, volatility estimation from high-frequen\-cy data has
attracted a lot of attention in financial econometrics and statistics.
Due to empirical evidence that the observed transaction prices of
assets cannot follow a discretely sampled semi-martingale model, a
prominent approach is to model the observations as the superposition of
the true (or efficient) price process with some measurement error,
conceived as microstructure noise. Main features are already present in
the basic model of observing
%
%e1.1 ###
\begin{equation}\label{EqE1Obs}
Y_i=X_{i/n}+\eps_i,\qquad   i=1,\ldots,n,
\end{equation}
with an efficient price process $X_t=\int_0^t\sigma(s) \,dB_s$, $B$ a
standard Brownian motion, and $\eps_i\sim N(0,\delta^2)$ all
independent. The aim is to perform
statistical inference on the volatility function $\sigma\dvtx[0,1]\to\R
^+$, for example, estimating the so-called integrated volatility $\int
_0^1\sigma^2(t) \,dt$ over the trading day.

The mathematical foundation on the parametric formulation of this model
has been laid by \cit{GloterJacodI} who prove the interesting result
that the model is locally asymptotically normal (LAN) as $n\to\infty
$, but with the unusual rate $n^{-1/4}$, while without microstructure
noise the rate is $n^{-1/2}$. Starting with \cit{Zhangetal}, the
nonparametric model has come into the focus of research. Mainly three
different, but closely related approaches have been proposed afterwards
to estimate the integrated volatility: multi-scale estimators [\citet{Zhang}],
realized kernels or autocovariances [\citet{BNetal}] and
preaveraging [\citet{Jacodetal}].
Under various degrees of generality, especially also for stochastic
volatility, all authors provide central limit theorems with convergence
rate $n^{-1/4}$ and an asymptotic variance involving the so-called
quarticity $\int_0^1\sigma^4(t) \,dt$. Recently, also rate-optimal
estimators for the spot volatility $\sigma^2(t)$ have been proposed
[\citet{MunkJohannes}, \citet{HoffMunkJohannes}].

The aim of the present paper is to provide a thorough mathematical
understanding of the basic model, to explain more profoundly why
statistical inference is not so canonical and to propose a simple
estimator of the integrated volatility which is efficient. To this end,
we employ Le Cam's concept of asymptotic equivalence between
experiments. In fact, our main theoretical result in Theorem \ref
{ThmE0G0} states under the $\alpha$-H\"{o}lder-regularity condition
$\alpha\ge(1+\sqrt{5})/4\approx0.81$ for $\sigma^2(\bullet)$ that
observing $(Y_i)$ in \eqref{EqE1Obs} is for $n\to\infty$
asymptotically equivalent to observing the Gaussian shift experiment
\[
dY_t=\sqrt{2\sigma(t)}\, dt+\delta^{1/2}n^{-1/4} \,dW_t, \qquad  t\in[0,1],
\]
with Gaussian white noise $dW$.
%Abstractly, this means that we observe for all functions $\phi\in
%L^2([0,1])$ the random variables $\int\phi dY=\int\phi(t)\sqrt{2
%N(0,1)$ and $\Cov(\eta_\phi,\eta_\psi)=\scapro{\phi}{\psi}_{L^2}$, see
%also Section \ref{SecPrelim} below.
By the \cit{BrownLow} result, we~ob\-tain a fortiori asymptotic
equivalence with the regression model
\[
Y_i=\sqrt{2\sigma\bigl(i/\sqrt{n}\bigr)}+\delta^{1/2}\eps_i,\qquad   i=1,\ldots
,\sqrt{n},\ \eps_i\sim N(0,1)\mbox{ i.i.d.}
\]

Not only the large noise level $\delta^{1/2}n^{-1/4}$ is apparent, but
also a nonlinear $\sqrt{\sigma(t)}$-form of the signal, from which
optimal asymptotic variance results can be derived. Note that a similar
form of a Gaussian shift was found to be asymptotically equivalent to
nonparametric density estimation [\citet{Nussbaum}]. A key ingredient of
our asymptotic equivalence proof are the results by \cit
{GramaNussbaum02} on asymptotic equivalence for generalized
nonparametric regression, but also ideas from \cit{Carter} and \cit
{Reiss08} play a role. Moreover, fine bounds on Hellinger distances for
Gaussian measures with different covariance operators turn out to be essential.

Roughly speaking, asymptotic equivalence means that any statistical
inference procedure can be transferred from one experiment to the other
such that the asymptotic risk remains the same, at least for bounded
loss functions. Technically, two sequences of experiments $\mathcal{E}^n$
and $\mathcal{G}^n$, defined on possibly different sample spaces, but with
the same parameter set, are asymptotically equivalent if the Le Cam
distance $\Delta(\mathcal{E}^n,\mathcal{G}^n)$ tends to zero. For $\mathcal{E}_i=(\mathcal{X}_i,\mathcal{F}_i,(\PP_\theta^i)_{\theta\in\Theta})$,
$i=1,2$, by definition, $\Delta(\mathcal{E}_1,\mathcal{E}_2)=\max(\delta
(\mathcal{E}_1,\mathcal{E}_2),\delta(\mathcal{E}_2,\mathcal{E}_1))$ holds in terms
of the deficiency $\delta(\mathcal{E}_1,\mathcal{E}_2)=\inf_M\sup_{\theta
\in\Theta}\norm{M\PP_\theta^1-\PP_\theta^2}_{\mathrm{TV}}$, where the
infimum is taken over all randomisations or Markov kernels $M$ from
$(\mathcal{X}_1,\mathcal{F}_1)$ to $(\mathcal{X}_2,\mathcal{F}_2)$; see, for
example, \cit{LeCamYang} for details. In particular, $\delta(\mathcal{E}_1,\mathcal{E}_2)=0$ means that $\mathcal{E}_1$ is more informative than
$\mathcal{E}_2$ in the sense that any observation in $\mathcal{E}_2$ can be
obtained from $\mathcal{E}_1$, possibly using additional randomizations.
Here, we shall always explicitly construct the transformations and
randomizations and we shall then only use that $\Delta(\mathcal{E}_1,\mathcal{E}_2)\le\sup_{\theta\in\Theta}\norm{\PP_\theta^1-\PP
_\theta^2}_{\mathrm{TV}}$ holds when both experiments are defined on the same
sample space.

The asymptotic equivalence is deduced stepwise. In Section \ref
{SecRegWNM}, the regres\-sion-type model \eqref{EqE1Obs} is shown to be
asymptotically equivalent to a corres\-ponding white noise model with
signal $X$. Then in Section \ref{SecLessInfG}, a very \mbox{simple}
construction yields a Gaussian shift model with signal $\log(\sigma
^2(\bullet)+c)$, $c>0$ some constant, which is asymptotically less
informative, but only by a constant factor in the Fisher information.
Inspired by this construction, we~pre\-sent a generalization in Section
\ref{SecSeqSimple} where the information loss can be made~ar\-bitrarily
small (but not zero), before applying nonparametric local asymptotic
theory in Section \ref{SecLoc} to derive asymptotic equivalence with
our final Gaussian shift model for shrinking local neighborhoods of
the parameters. Section \ref{SecGlob} yields the global result, which
is based on an asymptotic sufficiency result for simple independent statistics.

Extensions and restrictions are discussed in Section \ref{SecDisc},
where we also~pre\-sent a counter-example which shows that asymptotic
equivalence fails for~H\"{o}l\-der smoothness $\alpha< 1/3$ of the
volatility function $\sigma^2(\bullet)$. To determine whether
\mbox{asymptotic} equivalence holds or fails for $\alpha\in[1/3,(1+\sqrt
{5})/4]$ remains a~chal\-lenging open problem. In Section \ref{SecAppl},
we use the theoretical insight to construct a rate-optimal estimator of
the spot volatility and an efficient estimator of the integrated
volatility by a genuine local-likelihood approach. Remarkably, the
asymptotic variance is found to depend on the third moment $\int
_0^1\sigma^3(t) \,dt$ and for nonconstant $\sigma^2(\bullet)$ our
estimator outperforms previous approaches applied to the basic model.
Constructions needed for the proof are presented and discussed
alongside the mathematical results, deferring more technical parts to
the \hyperref[app]{Appendix}, which in Section~\ref{SecPrelim} also contains a summary
of results on white noise models, the Hellinger distance and
Hilbert--Schmidt norm estimates.

%s2 ###
\section{The regression and white noise model}\label{SecRegWNM}

In the main part, we shall work in the white noise setting, which is
more intuitive to handle than the regression setting, which in turn is
the observation model in practice. Let us define both models formally.
For that, we introduce the H\"{o}lder ball
\begin{eqnarray}
C^\alpha(R):=\{ f\in C^\alpha([0,1]) | \norm{f}_{C^\alpha}\le R\}\nonumber
\\
\eqntext{\displaystyle \mbox{with }
\norm{f}_{C^\alpha}=\norm{f}_\infty+\sup_{x\not=y}\frac{\abs
{f(x)-f(y)}}{\abs{x-y}^\alpha}.}
\end{eqnarray}

\begin{definition}
Let $\mathcal{E}_0=\mathcal{E}_0(n,\delta,\alpha,R,\underline{\sigma}^2)$
with $n\in\N$, $\delta>0$, $\alpha\in(0,1)$, $R>0$, $\underline
{\sigma}^2\ge0$ be the statistical experiment generated by observing
\eqref{EqE1Obs}. The volatility $\sigma^2$ belongs to the class
\[
\mathcal{S}(\alpha,R,\underline{\sigma}^2):= \Bigl\{\sigma^2\in
C^\alpha(R)  \bigm|  \min_{t\in[0,1]}\sigma^2(t)\ge\underline
{\sigma}^2 \Bigr\}.
\]

Let $\mathcal{E}_1=\mathcal{E}_1(\eps,\alpha,R,\underline{\sigma}^2)$
with $\eps>0$, $\alpha\in(0,1)$, $R>0$, $\underline{\sigma}^2\ge
0$ be the statistical experiment generated by observing
\[
dY_t=X_t \,dt+\eps \,dW_t,  \qquad t\in[0,1],
\]
with $X_t=\int_0^t\sigma(s) \,dB_s$ as above, independent standard
Brownian motions $W$ and $B$ and $\sigma^2\in\mathcal{S}(\alpha
,R,\underline{\sigma}^2)$.
\end{definition}

From \cit{BrownLow}, it is well known that the white noise and the
Gaussian regression model are asymptotically equivalent for noise level
$\eps=\delta/\sqrt{n}\to0$ as $n\to\infty$, provided the signal
is $\beta$-H\"{o}lder continuous for $\beta>1/2$.
Since Brownian motion and thus also our underlying process $X$ is only
H\"{o}lder continuous of order $\beta<1/2$ (whatever $\alpha$ is), it
is not clear whether asymptotic equivalence can hold for the
experiments $\mathcal{E}_0$ and $\mathcal{E}_1$. Yet, this is true.
Subsequently, we employ the notation $A_n\lesssim B_n$ if $A_n=O(B_n)$
and $A_n\thicksim B_n$ if $A_n\lesssim B_n$ as well as $B_n\lesssim
A_n$ and obtain the following theorem.

\begin{theorem}\label{ThmRegression}
For any $\alpha>0$, $\underline\sigma^2\ge0$ and $\delta,R>0$ the
experiments $\mathcal{E}_0$ and $\mathcal{E}_1$ with $\eps=\delta/\sqrt
{n}$ are asymptotically equivalent; more precisely,
\[
\Delta\bigl(\mathcal{E}_0(n,\delta,\alpha,R,\underline{\sigma}^2),
\mathcal{E}_1\bigl(\delta/\sqrt{n},\alpha,R,\underline{\sigma}^2\bigr)\bigr)\lesssim
R\delta^{-2}n^{-\alpha}.
\]
\end{theorem}

Interestingly, the asymptotic equivalence holds for any positive H\"
{o}lder regularity $\alpha>0$. In particular, for this result the
volatility $\sigma^2$ could be itself a continuous semi-martingale,
but such that $X$ conditionally on $\sigma^2$ remains Gaussian. Let us
also recall that by inclusion asymptotic equivalence always holds for
subclasses of functions, here for example for $C^m$-balls of $m$-times
continuously differentiable functions $\sigma^2$ so that we write
$\alpha>0$, meaning arbitrarily small positive $\alpha$, and not
$\alpha\in(0,1]$, which is more formal, but misleading. As the proof
in Section \ref{SecProofThmRegression} of the \hyperref[app]{Appendix} reveals, we
construct the equivalence by rate-optimal approximations of the
anti-derivative of $\sigma^2$ which lies in $C^{1+\alpha}$. Similar
techniques have been used by \cit{Carter} and \cit{Reiss08}, but here
we have to cope with the random signal for which we need to bound the
Hilbert--Schmidt norm of the respective covariance operators. Note
further that the asymptotic equivalence even holds when the noise level
$\delta$ tends to zero, provided $\delta^2n^\alpha\to\infty$
remains valid.

%s3 ###
\section{Less informative Gaussian shift experiments}\label{SecLessInfG}

From now on, we shall work with the white noise observation experiment
$\mathcal{E}_1$, where the main structures are more clearly visible. In
this section, we shall find easy Gaussian shift models which are
asymptotically not more informative than $\mathcal{E}_1$, but already
permit rate-optimal estimation results. The whole idea is easy to grasp
once we can replace the volatility $\sigma^2$ by a piecewise constant
approximation on small blocks of size $h$. That this is no loss of
generality is shown by the subsequent asymptotic equivalence result,
proved in Section \ref{SecPropPiecewiseConstant} of the \hyperref[app]{Appendix}.

\begin{definition}
Let $\mathcal{E}_2=\mathcal{E}_2(\eps,h,\alpha,R,\underline{\sigma}^2)$
be the statistical experiment generated by observing
\[
dY_t=X_t^h dt+\eps \,dW_t, \qquad  t\in[0,1],
\]
with $X_t^h=\int_0^t\sigma(\floor{s}_h) \,dB_s$, $\floor
{s}_h:=\floor{s/h}h$ for $h>0$ and $h^{-1}\in\N$, and independent
standard Brownian motions $W$ and $B$. The volatility $\sigma^2$
belongs to the class $\mathcal{S}(\alpha,R,\underline{\sigma}^2)$.
\end{definition}

\begin{proposition}\label{PropPiecewiseConstant}
Assume $\alpha\in(1/2,1]$ and $\underline\sigma^2>0$. Then for
$\eps\to0$, $h^\alpha=o(\eps^{1/2})$ the experiments $\mathcal{E}_1$
and $\mathcal{E}_2$ are asymptotically equivalent; more precisely,\vspace{-5pt}
\[
\Delta\bigl(\mathcal{E}_1(\eps,\alpha,R,\underline{\sigma}^2),
\mathcal{E}_2(\eps,h,\alpha,R,\underline{\sigma}^2)\bigr)\lesssim R
\underline{\sigma}^{-3/2} h^\alpha\eps^{-1/2}.
\]
\end{proposition}

In the sequel, we always assume $h^\alpha=o(\eps^{1/2})$ to hold such
that we can work equivalently with $\mathcal{E}_2$. Recall that observing
$Y$ in a white noise model is equivalent to observing $(\int e_m
\,dY)_{m\ge1}$ for an orthonormal basis $(e_m)_{m\ge1}$ of
$L^2([0,1])$; cf. also Section \ref{SecPrelim} below. Our first
step is thus to find an orthonormal system (not a basis) which extracts
as much {\it local } information on $\sigma^2$ as possible. For any
$\phi\in L^2([0,1])$ with $\norm{\phi}_{L^2}=1$, we have by partial
integration
%
%e3.1 ###
\begin{eqnarray} \label{EqCovY}
\qquad \int_0^1 \phi(t) \,dY_t&=&\int_0^1 \phi(t)X_t^h \,dt+\eps\int_0^1\phi
(t) \,dW_t\nonumber\\
&=&\Phi(1)X_1^h-\Phi(0)X_0^h-\int_0^1 \Phi(t)\sigma(\floor{t}_h)
\,dB_t+\eps\int\phi(t) \,dW_t\\
&=& \biggl(\int_0^1\Phi^2(t)\sigma^2(\floor{t}_h)\, dt+\eps^2
\biggr)^{1/2}\zeta_\phi,\nonumber
\end{eqnarray}
where $\Phi(t)=-\int_t^1\phi(s) \,ds$ is the antiderivative of $\phi
$ with $\Phi(1)=0$ and $\zeta_\phi\sim N(0,1)$ holds. To ensure that
$\Phi$ has only support in some interval $[kh,(k+1)h]$, we require
$\phi$ to have support in $[kh,(k+1)h]$ and to satisfy $\int\phi(t)
\,dt=0$. The function $\phi_k$ with $\supp(\phi_k)=[kh,(k+1)h]$,
$\norm{\phi_k}_{L^2}=1$, $\int\phi_k(t) \,dt=0$ that maximizes the
information load $\int\Phi_k^2(t) \,dt$ for $\sigma^2(kh)$ is given
by (use Lagrange theory)
%
%e3.2 ###
\begin{equation}\label{Eqphik}
\phi_k(t)=\sqrt{2}h^{-1/2}\cos \bigl(\pi(t-kh)/h \bigr) \mathbf{
1}_{[kh,(k+1)h]}(t),\qquad   t\in[0,1].
\end{equation}
The $L^2$-orthonormal system $(\phi_k)$ for $k=0,1,\ldots,h^{-1}-1$
is now used to construct Gaussian shift observations. In $\mathcal{E}_2$,
we obtain from \eqref{EqCovY} the observations
%
%e3.3 ###
\begin{equation}\label{Eqyk}
\qquad y_k:=\int\phi_k(t)\, dY_t= \bigl(h^2\pi^{-2}\sigma^2(kh)+\eps^2
\bigr)^{1/2}\zeta_k, \qquad  k=0,\ldots,h^{-1}-1,
\end{equation}
with independent standard normal random variables $(\zeta
_k)_{k=0,\ldots,h^{-1}-1}$. Observing $(y_k)$ is equivalent to observing
%
%e3.4 ###
\begin{equation}\label{Eqzk}
z_k:=\log(y_k^2h^{-2}\pi^2)-\E[\log(\zeta_k^2)]
=\log \bigl(\sigma^2(kh)+\eps^2h^{-2}\pi^2 \bigr)+\eta_k
\end{equation}
for $k=0,\ldots,h^{-1}-1$ with $\eta_k:=\log(\zeta_k^2)-\E[\log
(\zeta_k^2)]$ since $(y_k^2)$ is a sufficient statistic in \eqref
{Eqyk} and the logarithm is one-to-one.

We have found a nonparametric regression model with regression function
$\log(\sigma^2(\bullet)+\eps^2h^{-2}\pi^2)$ and $h^{-1}$ equidistant
observations corrupted by non-Gaussian, but centered noise $(\eta_k)$
of variance 2. To ensure that the regression function does not change
under the asymptotics $\eps\to0$, we specify the block size $h=h(\eps
)=h_0\eps$ with some fixed constant $h_0>0$.

It is not surprising that the nonparametric regression experiment in
\eqref{Eqzk} is equivalent to a corresponding Gaussian shift
experiment. Indeed, this follows readily from results by \cit
{GramaNussbaum02} who in their Section 4.2 derive asymptotic
equivalence already for our Gaussian scale model \eqref{Eqyk}. Note,
however, that their Fisher information for $\theta=\sigma^2$ must be
corrected to $I(\theta)=\frac12 \theta^{-2}$. We then obtain
directly asymptotic equivalence of \eqref{Eqyk} with the Gaussian
regression model
\[
w_k=\frac{1}{\sqrt{2}}\log\bigl(\sigma^2(kh)+h_0^{-2}\pi^2\bigr)+\gamma
_k, \qquad  k=0,\ldots,h^{-1}-1,
\]
where $\gamma_k\sim N(0,1)$ i.i.d. Since by the classical result of
\cit{BrownLow} or by \cit{Reiss08} the Gaussian regression is
equivalent to the corresponding white noise experiment [note that $\log
(\sigma^2(\bullet)+h_0^{-2}\pi^2)$ is also $\alpha$-H\"{o}lder
continuous], we have already derived an important and far-reaching result.

\begin{theorem}\label{ThmEquivE1G1}
For $\alpha>1/2$ and $\underline\sigma^2>0$ the high-frequency
experiment $\mathcal{E}_1(\eps,\alpha,R,\underline{\sigma}^2)$ is
asymptotically more informative than the Gaussian shift experiment
$\mathcal{G}_1(\eps,\alpha,R,\underline{\sigma}^2,h_0)$ of observing
\[
dZ_t=\frac{1}{\sqrt{2}}\log \bigl(\sigma^2(t)+h_0^{-2}\pi^2 \bigr)\, dt
+h_0^{1/2}\eps^{1/2}\,dW_t,\qquad   t\in[0,1].
\]
Here $h_0>0$ is an arbitrary constant and $\sigma^2\in\mathcal{S}(\alpha
,R,\underline{\sigma}^2)$.
\end{theorem}

\begin{remark}\label{RemG1loss}
Moving the constants from the diffusion to the drift part, the
experiment $\mathcal{G}_1$ is equivalent to observing
%
%e3.5 ###
\begin{equation}\label{EqZtilde}
d\tilde Z_t=(2h_0)^{-1/2}\log\bigl(\sigma^2(t)+h_0^{-2}\pi^2\bigr)\, dt
+\eps^{1/2}\,dW_t,\qquad   t\in[0,1].
\end{equation}
%
%The Gaussian shift experiment is nonlinear in $\sigma^2$ which is to
%be expected.
Writing $\eps=\delta/\sqrt{n}$ gives us the noise level $\delta
^{1/2}n^{-1/4}$ which appears in all~pre\-vious work on the model $\mathcal{E}_0$.

To quantify the amount of information we have lost, let us study the
LAN-property of the constant parametric case $\sigma^2(t)=\sigma^2>0$
in $\mathcal{G}_1$. We consider the local alternatives $\sigma_\eps
^2=\sigma_0^2+\eps^{1/2}$ for which we obtain the Fisher information
$I_{h_0}=(2h_0)^{-1}h_0^4/(\pi^2+h_0^2\sigma_0^2)^2$. Maximizing over
$h_0$ yields $h_0=\sqrt{3}\pi\sigma_0^{-1}$ and the Fisher
information is at most equal to $\sup_{h_0>0}I_{h_0}=\sigma
_0^{-3}3^{3/2}/(32\pi)\approx0.0517\sigma_0^{-3}$.

By the LAN-result of \cit{GloterJacodI} for $\mathcal{E}_0$, the best
value~is $I(\sigma_0)=\frac18\sigma_0^{-3}$ which is clearly larger.
Note, however, that the relative (normalized) efficiency is already
$\frac{\sqrt{3^{3/2}/(32\pi)}}{\sqrt{1/8}}\approx0.64$, which
means that we attain here about $64\%$ of the precision when working
with $\mathcal{G}_1$ instead of $\mathcal{E}_0$ or $\mathcal{E}_1$.
\end{remark}

%s4 ###
\section{A close sequence of simple models}\label{SecSeqSimple}

In order to decrease the information loss in $\mathcal{G}_1$, we now take
into account higher frequencies in each block $[kh,(k+1)h]$ by using
further trigonometric basis functions. In the case of constant $\sigma
^2$, the covariance operator of the observations is diagonalized by the
Karhunen--Lo\`{e}ve basis for Brownian motion which together with a~block\-wise
approximation is exactly the idea here; see also the
discussion in Section~\ref{SecDisc}. Equivalently, we can argue by a
variational principle, maximizing the information load as in the case
of $\phi_k$.
%We also allow at this point for general $h=o(\eps^\alpha)$, no longer
%requiring the %order $\eps$.
In a frequency-location notation $(j,k)$, we consider for $k=0,1,\ldots
,h^{-1}-1, j\ge1$,
%
%e4.1 ###
\begin{equation}\label{Eqphijk}
\qquad \phi_{jk}(t)=\sqrt{2}h^{-1/2}\cos\bigl(j\pi(t-kh)/h\bigr)\mathbf{1}_{[kh,(k+1)h]}(t),\qquad   t\in[0,1].
\end{equation}
This gives the corresponding antiderivatives
\[
\Phi_{jk}(t)=\frac{\sqrt{2h}}{\pi j}\sin\bigl(j\pi(t-kh)/h\bigr)\mathbf{1}_{[kh,(k+1)h]}(t),\qquad   t\in[0,1].
\]
Not only the $(\phi_{jk})$ and $(\Phi_{jk})$ are localized on each
block, also each single~fa\-mily of functions is orthogonal in
$L^2([0,1])$. Working again on the piecewise constant experiment $\mathcal{E}_2$, we extract the observations
%
%e4.2 ###
\begin{eqnarray} \label{Eqyjk}
y_{jk}:=\int_0^1 \phi_{jk}(t)\, dY_t= \bigl(h^2\pi^{-2}j^{-2}\sigma
^2(kh)+\eps^2 \bigr)^{1/2}\zeta_{jk},\nonumber\\[-8pt]\\[-8pt]
\eqntext{j\ge1, k=0,\ldots,h^{-1}-1,}
\end{eqnarray}
with $\zeta_{jk}\sim N(0,1)$ independent over all $(j,k)$. Note that
independence follows since $(\phi_{jk})$ and $(\Phi_{jk})$ are both
$L^2$-orthogonal families and the observations are therefore
uncorrelated. The same transformation as before leads for each $j\ge1$
to the regression model for $k=0,\ldots,h^{-1}-1$
%
%e4.3 ###
\begin{eqnarray} \label{Eqzjk}
z_{jk}&:=&\log(y_{jk}^2)-\log(h^2\pi
^{-2}j^{-2})-\E[\log(\zeta_{jk}^2)]\nonumber\\[-8pt]\\[-8pt]
&\hspace*{3pt}=&\log\bigl(\sigma^2(t)+\eps^2h^{-2}\pi^2j^2\bigr)+\eta_{jk}.\nonumber
\end{eqnarray}
Applying the asymptotic equivalence result by \cit{GramaNussbaum02}
for each independent level $j$ separately, we immediately generalize
Theorem~\ref{ThmEquivE1G1}.

\begin{theorem}\label{ThmE1G2}
For $\alpha>1/2$ and $\underline{\sigma}^2>0$, the high-frequency
experiment $\mathcal{E}_1(\eps,\alpha,R,\underline{\sigma}^2)$ is
asymptotically more informative than the combined experiment $\mathcal{G}_2(\eps,\alpha,R,\underline{\sigma}^2,h_0,J)$ of independent
Gaussian shifts
\begin{eqnarray}
dZ^j_t=\frac{1}{\sqrt{2}}\log\bigl(\sigma^2(t)+h_0^{-2}\pi^2j^2\bigr)\, dt
+h_0^{1/2}\eps^{1/2}\,dW^j_t, \nonumber \\
\eqntext{t\in[0,1], j=1,\ldots,J,}
\end{eqnarray}
with independent Brownian motions $(W^j)_{j=1,\ldots,J}$ and $\sigma
^2\in\mathcal{S}(\alpha,R,\underline{\sigma}^2)$.
%, which in turn is equivalent to observing
% dt
%+\eps^{1/2}dW^j_t,  t\in[0,1], j\ge1.
The constants $h_0>0$ and $J\in\N$ are arbitrary, but fixed.
\end{theorem}

\begin{remark}
Let us again study the LAN-property of the constant~pa\-rametric case
$\sigma^2(t)=\sigma^2>0$ for the local alternatives $\sigma^2_\eps
=\sigma_0^2+\eps^{1/2}$. We obtain the Fisher information
\[
I_{h_0,J}=\sum_{j=1}^J(2h_0)^{-1}h_0^4(\pi^2j^2+h_0^2\sigma_0^2)^{-2}
=\sum_{j=1}^J\frac{h_0^{-1}}{2(\pi^2(jh_0^{-1})^2+\sigma_0^2)^2}.
\]
In the limit $J\to\infty$ and $h_0\to\infty$, we obtain by Riemann
sum approximation
\[
\lim_{h_0\to\infty}\lim_{J\to\infty} I_{h_0,J}=\int_0^\infty
\frac{dx}{2(\pi^2x^2+\sigma_0^2)^2}
%= (\frac{x}{4\sigma_0^2(\sigma_0^2+\pi^2x^2)}+\frac{\arctan(\pi x/
=\frac{1}{8\sigma_0^3}.
\]
This is exactly the optimal Fisher information, obtained by \cit
{GloterJacodI} in this case. Note, however, that it is not at all
obvious that we may let $J,h_0\to\infty$, in the asymptotic
equivalence result. Moreover, in our theory the restriction $h^\alpha
=o(\eps^{1/2})$ is necessary, which translates into $h_0=o(\eps
^{(1-2\alpha)/2\alpha})$. Still, the positive aspect is that we can
come as close as we wish to an asymptotically almost equivalent, but
much simpler model. The convergence $h_0\to\infty$ is also an
essential point in the final proof, starting with the next section.
\end{remark}

%s5 ###
\section{Localization}\label{SecLoc}

We know from standard regression theory [\citet{Stone82}] that in the
experiment $\mathcal{G}_1$ we can estimate $\sigma^2\in C^\alpha$ in
sup-norm with rate $(\eps\log(\eps^{-1}))^{\alpha/(2\alpha+1)}$,
using that the log-function is a $C^\infty$-diffeomorphism for
arguments bounded away from zero and infinity. Since $\mathcal{E}_1$ is
for $\alpha>1/2$ asymptotically more informative than $\mathcal{G}_1$, we
can therefore localize $\sigma^2$ in a~neighborhood of some $\sigma
_0^2$. Using the local coordinate $s^2$ in $\sigma^2=\sigma
_0^2+v_\eps s^2$ for $v_\eps\to0$, we define a localized experiment;
cf. \cit{Nussbaum}.

\begin{definition}
Let $\mathcal{E}_{i,\mathrm{loc}}=\mathcal{E}_{i,\mathrm{loc}}(\sigma_0,\eps,\alpha
,R,\underline{\sigma}^2)$ for $\sigma_0\in\mathcal{S}(\alpha
,R,\underline{\sigma}^2)$ be the statistical subexperiment obtained
from $\mathcal{E}_i(\eps,\alpha,R,\underline{\sigma}^2)$ by
restricting to the parameters $\sigma^2=\sigma_0^2+v_\eps s^2$ with
$v_\eps=\eps^{\alpha/(2\alpha+1)}\log(\eps^{-1})$ and unknown
$s^2\in C^\alpha(R)$.
\end{definition}

We shall consider the observations $(y_{jk})$ in \eqref{Eqyjk} derived
from $\mathcal{E}_{2,\mathrm{loc}}$ and multiplied by $\pi j/h$. The model is then
a generalized nonparametric regression family in the sense of \cit
{GramaNussbaum02}. On the sequence space $(\mathcal{X},\mathcal{F})=(\R^{\N
},{\mathfrak B}^{\otimes\N})$, we consider for $\theta\in\Theta
=[\underline\sigma^2,R]$ the Gaussian product measure
%
%e5.1 ###
\begin{equation}\label{EqPtheta}
\PP_\theta=\bigotimes_{j\ge1} N (0,\theta+h_0^{-2}\pi
^2j^2 ).
\end{equation}
The parameter $\theta$ plays the role of $\sigma^2(kh)$ for each $k$.
By independence and the result for the one-dimensional Gaussian scale
model, the Fisher information for $\theta$ is given by
%
%e5.2 ###
\begin{eqnarray} \label{EqItheta}
I(\theta)&:=&\sum_{j\ge1} \frac
{1}{2(\theta+h_0^{-2}\pi^2j^2)^2}\nonumber\\[-8pt]\\[-8pt]
&\hspace*{3pt}=&\frac{h_0}{8\theta^{3/2}} \biggl(\frac{1+4\theta^{1/2} h_0
e^{-2\theta^{1/2} h_0}-e^{-4\theta^{1/2} h_0}}{(1-e^{-2\theta^{1/2} h_0})^2}
-\frac{2}{\theta^{1/2} h_0} \biggr),\nonumber
\end{eqnarray}
where the series is evaluated using the derivative with respect to
$\alpha$ in the identity
$\sum_{j=1}^\infty\frac{1}{j^2+\alpha^2}=\frac{1}{2\alpha^2}(\pi
\alpha\coth(\pi\alpha)-1)$.
Since we shall later let $h_0$ tend to infinity, an essential point is
the asymptotics $I(\theta)\thicksim h_0$.

We split our observation design $\{kh \mid k=0,\ldots,h^{-1}\}$ into
blocks $A_m=\{kh \mid k=(m-1)\ell,\ldots,m\ell-1\}$, $m=1,\ldots
,(\ell h)^{-1}$, of length $\ell$ such that the radius $v_\eps$ of
our nonparametric local neighborhood has the order of the {\it
parametric} noise level $(I(\theta)\ell)^{-1/2}$ in each block:
%
%e5.3 ###
\begin{equation}\label{EqParRate}
v_\eps\thicksim(I(\theta)\ell)^{-1/2}\quad\Rightarrow\quad\ell\thicksim
h_0^{-1}v_\eps^{-2}.%=o (h_0^{-1}(\eps\log(\eps^{-1}))^{\alpha/(2
\end{equation}

For later convenience, we consider odd and even indices $k$ separately,
assuming that $h^{-1}$ and $\ell$ are even integers.
This way, for each block $m$~ob\-serving $(y_{jk}\pi j/h)$ for $j\ge1$
and $k\in A_m$, $k$ odd, respectively, $k$ even, can be modeled by the
experiments
\begin{eqnarray} \label{EqE3m}
\mathcal{E}_{3,m}^{\mathrm{odd}}&=& \biggl(\mathcal{X}^{\ell/2},\mathcal{F}^{\otimes\ell
/2}, \biggl(\bigotimes_{k\in A_m\ \mathrm{\mathrm{odd}}}\PP_{\sigma
_0^2(k/n)+v_\eps s^2(k/n)} \biggr)_{s^2\in C^\alpha(R)} \biggr),\\
\mathcal{E}_{3,m}^{\mathrm{even}}&=& \biggl(\mathcal{X}^{\ell/2},\mathcal{F}^{\otimes\ell
/2}, \biggl(\bigotimes_{k\in A_m\ \mathrm{\mathrm{even}}}\PP_{\sigma
_0^2(k/n)+v_\eps s^2(k/n)} \biggr)_{s^2\in C^\alpha(R)} \biggr),
\end{eqnarray}
where all parameters are the same as for $\mathcal{E}_{2,\mathrm{loc}}$. Using the
nonparametric local asymptotic theory developed by \cit
{GramaNussbaum02} and the independence of the experiments $(\mathcal{E}_{3,m}^{\mathrm{odd}})_m$
[resp., $(\mathcal{E}_{3,m}^{\mathrm{even}})_m$], we are able to
prove in Section \ref{SecProofPropE1G3loc} the following asymptotic
equivalence.

\begin{proposition}\label{PropE1G3loc}
Assume $\alpha>1/2$, $\underline\sigma^2>0$ and $h_0\thicksim\eps
^{-p}$ with $p\in(0,1-(2\alpha)^{-1})$ such that $(2h)^{-1}\in\N$.
Then observing $\{y_{j,2k+1} \mid j\ge1, k=0,\ldots,\break (2h)^{-1}-1\}$ in
experiment $\mathcal{E}_{2,\mathrm{loc}}$ is asymptotically equivalent to the
local~Gaus\-sian shift experiment $\mathcal{G}_{3,\mathrm{loc}}$ of observing
%
%e5.4 ###
\begin{eqnarray} \label{EqG3loc}
dY_t=\frac{1}{\sqrt{8}\sigma_0^{3/2}(t)} \biggl(1
-\frac{2}{\sigma_0(t)h_0} \biggr)^{1/2}v_\eps s^2(t)\, dt+ (2\eps
)^{1/2\,}dW_t,\nonumber\\[-8pt]\\[-8pt]
\eqntext{t\in[0,1],}
\end{eqnarray}
where the unknown $s^2$ and all parameters are the same as in $\mathcal{E}_{2,\mathrm{loc}}$.
The Le Cam distance tends to zero uniformly over the
center of localization $\sigma_0^2\in\mathcal{S}(\alpha,R,\underline
{\sigma}^2)$.

The same asymptotic equivalence result holds true for observing $\{
y_{j,2k} \mid j\ge1, k=0,\ldots,(2h)^{-1}-1\}$ in experiment $\mathcal{E}_{2,\mathrm{loc}}$.
\end{proposition}

Note that in this model, combining even and odd indices $k$, we can
already infer the LAN-result by \cit{GloterJacodI}, but we still face
a second-order term of order $h_0^{-1}v_\eps$ in the drift. This term
is asymptotically negligible only if it is of smaller order than the
noise level $\eps^{1/2}$. To be able to choose $h_0$ sufficiently
large, we have to require a larger H\"{o}lder smoothness of the volatility.

\begin{corollary}\label{CorE1G4loc}
Assume $\alpha>\frac{1+\sqrt{17}}{8}\approx0.64$, $\underline
\sigma^2>0$ and $h_0\thicksim\eps^{-p}$ with $p\in(0,1-(2\alpha
)^{-1})$ such that $(2h)^{-1}\in\N$. Then observing $\{y_{j,2k+1} \mid
j\ge1, k=0,\ldots,(2h)^{-1}-1\}$ in experiment $\mathcal{E}_{2,\mathrm{loc}}$ is
asymptotically equivalent to the local Gaussian shift experiment $\mathcal{G}_{4,\mathrm{loc}}$ of observing
%
%e5.5 ###
\begin{equation}\label{EqG4loc}
dY_t=\frac{1}{\sqrt{8}\sigma_0^{3/2}(t)}v_\eps s^2(t)\, dt+ (2\eps
)^{1/2}\,dW_t, \qquad  t\in[0,1],
\end{equation}
where the unknown $s^2$ and all parameters are the same as in $\mathcal{E}_{2,\mathrm{loc}}$.
The Le Cam distance tends to zero uniformly over the
center of localization $\sigma_0^2\in\mathcal{S}(\alpha,R,\underline
{\sigma}^2)$.

The same asymptotic equivalence result holds true for observing $\{
y_{j,2k} \mid j\ge1, k=0,\ldots,(2h)^{-1}-1\}$ in experiment $\mathcal{E}_{2,\mathrm{loc}}$.
\end{corollary}

\begin{pf}
For $\alpha>\frac{1+\sqrt{17}}{8}$, the choice of $h_0=\eps^{-p}$
for some $p\in(\frac{1}{4\alpha+2},\frac{2\alpha-1}{2\alpha})$ is
possible and ensures that $h^{\alpha}=o(\eps^{1/2})$ holds as well as
$h_0^{-2}=o(v_\eps^{-2}\eps)$. Therefore, the Kullback--Leibler
divergence between the observations in $\mathcal{G}_3^{\mathrm{loc}}$ and in
$\mathcal{G}_4^{\mathrm{loc}}$ evaluates by the Cameron--Martin (or Girsanov) formula to
\[
\eps^{-1}\int_0^1 \frac{1}{8\sigma_0^{3}(t)}\biggl ( \biggl(1
-\frac{2}{\sigma_0(t)h_0} \biggr)^{1/2}-1 \biggr)^2v_\eps^2 s^4(t)
\,dt\lesssim\eps^{-1}h_0^{-2}v_\eps^2.
\]
Consequently, the Kullback--Leibler and thus also the total variation
distance tend to zero.
\end{pf}

In a last step, we find local experiments $\mathcal{G}_{5,\mathrm{loc}}$, which are
asymptotically equivalent to $\mathcal{G}_{4,\mathrm{loc}}$ and do not depend on
the center of localization $\sigma_0^2$. To this end, we use a
variance-stabilizing transform, based on the Taylor expansion
\[
\sqrt{2}x^{1/4}=\sqrt{2}x_0^{1/4}+\frac{1}{\sqrt
{8}}x_0^{-3/4}(x-x_0)+O\bigl((x-x_0)^2\bigr)
\]
which holds uniformly over $x,x_0$ on any compact subset of $(0,\infty
)$. Inserting $x=\sigma^2(t)=\sigma_0^2(t)+v_\eps s^2(t)$ and
$x_0=\sigma_0^2$ from our local model, we obtain
%
%e5.6 ###
\begin{equation}\label{EqIthetaTaylor}
\sqrt{2\sigma(t)}=\sqrt{2\sigma_0(t)}+\frac{1}{\sqrt{8}}\sigma
_0^{-3/2}(t)v_\eps s^2(t)+O(v_\eps^2).
\end{equation}

Since $v_\eps^2=o(\eps^{1/2})$ holds for $\alpha>1/2$, we can add
the uninformative signal $\sqrt{2}\sigma_0^{1/2}(t)$ to $Y$ in $\mathcal{G}_{4,\mathrm{loc}}$, replace the drift by
$\sqrt{2}\sigma^{1/2}(t)$ and still keep convergence of the total
variation distance, compare the preceding proof. Consequently, from
Corollary \ref{CorE1G4loc} we obtain the following result.

\begin{corollary}\label{CorE1G5loc}
Assume $\alpha>\frac{1+\sqrt{17}}{8}\approx0.64$, $\underline
\sigma^2>0$ and $h_0\thicksim\eps^{-p}$ with $p\in(0,1-(2\alpha
)^{-1})$ such that $(2h)^{-1}\in\N$. Then observing $\{y_{j,2k+1} \mid
j\ge1, k=0,\ldots,(2h)^{-1}-1\}$ in the experiment $\mathcal{E}_{2,\mathrm{loc}}$
is asymptotically equivalent to the local Gaussian shift experiment
$\mathcal{G}_{5,\mathrm{loc}}$ of observing
%
%e5.7 ###
\begin{equation}\label{EqG5loc}
dY_t=\sqrt{2\sigma(t)}\, dt+ (2\eps)^{1/2}\, dW_t, \qquad  t\in[0,1],
\end{equation}
where the unknown is $\sigma^2=\sigma_0^2+v_\eps s^2$ and all
parameters are the same as~in $\mathcal{E}_{2,\mathrm{loc}}$. The Le Cam distance
tends to zero uniformly over the center of locali\-zation $\sigma_0^2\in
\mathcal{S}(\alpha,R,\underline{\sigma}^2)$.

The same asymptotic equivalence result holds true for observing $\{
y_{j,2k} \mid j\ge1, k=0,\ldots,(2h)^{-1}-1\}$ in experiment $\mathcal{E}_{2,\mathrm{loc}}$.
\end{corollary}

%s6 ###
\section{Globalization}\label{SecGlob}

The globalization now basically follows the usual route, first
established by \cit{Nussbaum}. Essential for us is to show that
observing $(y_{jk})$ for $j\ge1$ is asymptotically sufficient in
$\mathcal{E}_2$. Then we can split the white noise observation experiment
$\mathcal{E}_2$ into two independent sub-experiments obtained from
$(y_{jk})$ for $k$ odd and $k$ even, respectively. Usually, a white
noise experiment can be split into two independent subexperiments with
the same drift and an increase by $\sqrt{2}$ in the noise level. Here,
however, this does not work since the two diffusions in the {\it
random} drift remain the same and thus independence fails.

Let us introduce the $L^2$-normalized step functions
\begin{eqnarray*}
\phi_{0,k}(t) &:=&(2h)^{-1/2} \bigl(\mathbf{1}_{[(k-1)h,kh]}(t)-\mathbf{1}_{[kh,(k+1)h]}(t) \bigr),\qquad   k=1,\ldots,h^{-1}-1,\\
\phi_{0,0}(t) &:=&h^{-1/2}\mathbf{1}_{[0,h]}(t).
\end{eqnarray*}
We obtain a normalized complete basis $(\phi_{jk})_{j\ge0,0\le k\le
h^{-1}-1}$ of $L^2([0,1])$ such that observing $Y$ in experiment $\mathcal{E}_2$ is equivalent to observing
\[
y_{jk}:=\int_0^1\phi_{jk}(t)\, dY_t, \qquad  j\ge0, k=0,\ldots,h^{-1}-1.
\]
Calculating the Fourier series, we can express the tent function $\Phi
_{0,k}$ with $\Phi_{0,k}'=\phi_{0,k}$ and $\Phi_{0,k}(1)=0$ as an
$L^2$-convergent series over the dilated sine functions $\Phi_{jk}$
and $\Phi_{j,k-1}$, $j\ge1$:
%
%e6.1 ###
\begin{equation}\label{EqPhi0k}
\quad \Phi_{0,k}(t)=\sum_{j\ge1}(-1)^{j+1}\Phi_{j,k-1}(t)+\sum_{j\ge
1}\Phi_{jk}(t),\qquad   k=1,\ldots, h^{-1}-1.\hspace*{-10pt}
\end{equation}
We also have $\Phi_{0,0}(t)=2\sum_{j\ge1}\Phi_{j,0}(t)$.
By partial integration, this implies (with $L^2$-convergence)
\begin{eqnarray}
\beta_{0,k} :=\scapro{\phi_{0,k}}{X}
=-\int_0^1 \Phi_{0,k}(t) \,dX(t)
=\sum_{j\ge1}(-1)^{j+1}\beta_{j,k-1}+\sum_{j\ge1}\beta_{jk}\nonumber \\
\eqntext{\mbox{where }\beta_{jk}:=\scapro{\phi_{jk}}{X}}
\end{eqnarray}
for $k\ge1$ and similarly $\beta_{0,0}=2\sum_{j\ge1}\beta_{j,0}$.
This means that the signal $\beta_{0,k}$ in $y_{0,k}$ can be perfectly
reconstructed from the signals in the $y_{j,k-1}$, $y_{jk}$.
For jointly Gaussian random variables, we obtain the conditional law in
$\mathcal{E}_2$
\[
\mathcal{L}(\beta_{jk} | y_{jk})
=N \biggl(\frac{\Var(\beta_{jk})}{\Var(y_{jk})}y_{jk},
\frac{\eps^2\Var(\beta_{jk})}{\Var(y_{jk})} \biggr),
%=N (\frac{\norm{\Phi_{jk}}^2\sigma^2(kh)}
%{\norm{\Phi_{jk}}^2\sigma^2(kh)+\eps^2}y_{jk}, \frac{\norm{
%{\norm{\Phi_{jk}}^2\sigma^2(kh)+\eps^2}\eps^2 ).
\]
which depends on the unknown $\sigma^2(kh)$. Given the results by \cit
{Stone82} and our less-informative Gaussian shift experiment $\mathcal{G}_1$ for $\alpha>1/2$, $\underline{\sigma}^2>0$, there is an
estimator $\hat\sigma^2_\eps$ based on $(y_{1,k})_k$ in $\mathcal{E}_2$ with
%
%e6.2 ###
\begin{equation}\label{EqConsistentEst}
\lim_{\eps\to0}\inf_{\sigma^2\in\mathcal{S}} \PP_{\sigma^2,\eps
}(\norm{\hat\sigma^2_\eps-\sigma^2}_\infty\le Rv_\eps)=1,
\end{equation}
where $v_\eps=\eps^{\alpha/(2\alpha+1)}\log(\eps^{-1})$ as in the
definitions of the localized experiments.

In a randomization step, we can thus generate independent
$N(0,1)$-distri\-buted random variables $\rho_{jk}$ to construct from
$(y_{jk})_{j\ge1,k}$
\[
\tilde\beta_{jk}:=\frac{\Var_\eps(\beta_{jk})}{\Var_\eps(y_{jk})}y_{jk}
+\frac{\eps\Var_\eps(\beta_{jk})^{1/2}}{\Var_\eps
(y_{jk})^{1/2}}\rho_{jk},\qquad   j\ge1,
\]
where the variance $\Var_\eps$ is the expression for $\Var$ where
the unknown values $\sigma^2(kh)$ are replaced by the estimated values
$\hat\sigma^2_\eps(kh)$:
%
%e6.3 ###
\begin{equation}\label{EqVareps}
\Var_\eps(y_{jk})=\Var_\eps(\beta_{jk})+\eps^2,\qquad  \Var_\eps
(\beta_{jk})=h^2\pi^{-2}j^{-2}\hat\sigma_\eps^2(kh).
\end{equation}
From this, we define $\tilde\beta_{0,k}:=\sum_{j\ge1}
((-1)^{j+1}\tilde\beta_{j,k-1}+\tilde\beta_{jk})$, $\tilde\beta
_{0,0}:=2\sum_{j\ge1}\tilde\beta_{j,0}$ and generate artificial
observations $(\tilde y_{0,k})$ such that the conditional law
$\mathcal{L}((\tilde y_{0,k})_k | (y_{jk})_{j\ge1,k})$ corresponds to
$\mathcal{L}(( y_{0,k})_k | (y_{jk})_{j\ge1,k})$ in the sense that it
is multivariate normal with mean $(\tilde\beta_{0k})_k$ and
(tri-diagonal) covariance matrix $\eps^2(\scapro{\phi_{0,k}}{\phi
_{0,k'}})_{k,k'}$.

In Section \ref{ProofPropE1locyjk}, we shall prove that the Hellinger
distance between the families of centered Gaussian random variables
$\mathcal{Y}:=\{y_{jk} \mid j\ge0, k=0,\ldots,h^{-1}-1\}$ and $\tilde
\mathcal{Y}:=\{\tilde y_{0,k} \mid k=0,\ldots,h^{-1}-1\}\cup\{y_{jk}
\mid
j\ge1, k=0,\ldots,h^{-1}-1\}$
tends to zero, provided $h_0^{-1}v_\eps^2=o(\eps)$, which is possible
when $\alpha>\frac{1+\sqrt{5}}{4}$ with the choice $h_0=\eps^{-p}$
for some $p\in(\frac{1}{2\alpha+1},\frac{2\alpha-1}{2\alpha})$.
In particular, this means that $(y_{jk})_{j\ge1,k}$ is asymptotically
sufficient and the information in $(y_{0,k})_k$ is asymptotically negligible.

\begin{proposition}\label{PropE1locyjk}
Assume $\alpha>\frac{1+\sqrt{5}}{4}\approx0.81$, $\underline\sigma
^2>0$ and $h^{-1}$ an even integer. Then the experiment $\mathcal{E}_2$ is
asymptotically equivalent to the product experiment $\mathcal{E}_{2,\mathrm{odd}}\otimes\mathcal{E}_{2,\mathrm{even}}$
 where $\mathcal{E}_{2,\mathrm{odd}}$ is
obtained from the observations $\{y_{j,2k+1} \mid j\ge1, k=0,\ldots
,(2h)^{-1}-1\}$ and $\mathcal{E}_{2,\mathrm{even}}$ from the observations $\{
y_{j,2k} \mid j\ge1, k=0,\ldots,(2h)^{-1}-1\}$ in experiment $\mathcal{E}_2$.
\end{proposition}

This key result permits to globalize the local result. In the sequel,
we always assume $\alpha>\frac{1+\sqrt{5}}{4}$ and $\underline
{\sigma}^2>0$. We start with the asymptotic equivalence between $\mathcal{E}_2$ and $\mathcal{E}_{2,\mathrm{odd}}\otimes\mathcal{E}_{2,\mathrm{even}}$. Using again an
estimator $\hat\sigma^2_\eps$ in $\mathcal{E}_{2,\mathrm{odd}}$ satisfying
\eqref{EqConsistentEst}, we can localize the second factor $\mathcal{E}_{2,\mathrm{even}}$ around $\hat\sigma_\eps^2$ and therefore by Corollary
\ref{CorE1G5loc} replace it by experiment $\mathcal{G}_{5,\mathrm{loc}}$; see
Theorem 3.2 in \cit{Nussbaum} for a formal proof. Since $\mathcal{G}_{5,\mathrm{loc}}$ does not depend on the center $\hat\sigma_\eps^2$, we
conclude that $\mathcal{E}_2$ is asymptotically equivalent to the product
experiment $\mathcal{E}_{2,\mathrm{odd}}\otimes\mathcal{G}_5$ where $\mathcal{G}_5$ has
the same parameters as $\mathcal{E}_2$ and is given by observing $Y$ in
\eqref{EqG5loc}. Now we use an estimator $\hat\sigma^2_\eps$ in
$\mathcal{G}_5$ satisfying \eqref{EqConsistentEst}, whose existence is
ensured by \cit{Stone82}, to localize $\mathcal{E}_{2,\mathrm{odd}}$. Corollary
\ref{CorE1G5loc} then allows again to replace the localized $\mathcal{E}_{2,\mathrm{odd}}$-experiment by $\mathcal{G}_5$ such that $\mathcal{E}_2$ is
asymptotically equivalent to the product experiment $\mathcal{G}_5\otimes
\mathcal{G}_5$. Finally, taking the mean of the independent observations
\eqref{EqG5loc} in both factors, which is a sufficient statistic (or,
abstractly, due to identical likelihood processes) we see that $\mathcal{G}_5\otimes\mathcal{G}_5$ is equivalent to the experiment $\mathcal{G}_0$ of
observing
$dY_t=\sqrt{2\sigma(t)} \,dt+ \sqrt{\eps}\, dW_t$, $t\in[0,1]$.
Our final result then follows from the asymptotic equivalence between
$\mathcal{E}_0$ and $\mathcal{E}_1$ as well as between $\mathcal{E}_1$ and
$\mathcal{E}_2$.

\begin{theorem}\label{ThmE0G0}
Assume $\alpha>\frac{1+\sqrt{5}}{4}\approx0.81$ and $\delta
_n$, $\underline\sigma^2$, $R>0$. Then the regression experiment $\mathcal{E}_0(n,\delta_n,\alpha,R,\underline{\sigma}^2)$ is for $n\to\infty
$ and $\delta_n^{-2}n^{-\alpha}\to0$ asymptotically equivalent to
the Gaussian shift experiment $\mathcal{G}_0(\delta n^{-1/2},\alpha
,R,\underline{\sigma}^2)$ of observing
%
%e6.4 ###
\begin{equation}\label{EqG0}
dY_t=\sqrt{2\sigma(t)}\, dt+ \delta^{1/2}n^{-1/4} \,dW_t,\qquad   t\in[0,1],
\end{equation}
for $\sigma^2\in\mathcal{S}(\alpha,R,\underline{\sigma}^2)$.
\end{theorem}

\section{Discussion}\label{SecDisc}

Our results show that inference for the volatility in the
high-frequency observation model under microstructure noise $\mathcal{E}_0$ is asymptotically as difficult as in the well-understood Gaussian
shift model $\mathcal{G}_0$. Remark that the constructions in
\citeauthor{GloterJacodI} (\citeyear{GloterJacodI,GloterJacodII}) rely on preliminary estimators at
the boundary of suitable blocks, while we require $\supp\Phi
_{jk}=[kh,(k+1)h]$ to obtain independence among blocks. In this
context, Proposition \ref{PropE1locyjk} shows asymptotic sufficiency
of observing only the increment process $X_t-X_{kh}$, $t\in
[kh,(k+1)h]$, on each block due to $\int\phi_{jk}(t) \,dt=0$ for $j\ge
1$. Naturally, the $(\phi_{jk})_{j\ge1}$ form exactly the
eigenfunctions of the covariance operator of Brownian motion on
$[kh,(k+1)h]$ and it suffices to use the block-wise Karhunen--Lo\`{e}ve
expansion for inference.

It should be remarked that a fortiori asymptotic equivalence also holds
when using instead of the ($\phi_{jk})$ different basis functions on
each block spanning the orthogonal complement of the constant functions
(i.e., integrating to zero). For practical applications, especially
when estimating the spot volatility curve, the blocking might produce
artifacts and wavelet bases which realize a well localized
time frequency analysis seem to be well suited, compare \cit{HoffMunkJohannes}.

It is interesting to note that both, model $\mathcal{E}_0$ and model
$\mathcal{G}_0$, are homogeneous in the sense that factors from the noise
(i.e., the $dW_t$-term) can be moved to the drift term and vice versa
such that, for example, high volatility can counterbalance a high noise
level $\delta$ or a large observation distance $1/n$. Another
phenomenon is that observing $\mathcal{E}_0$ $m$-times independently with
$n$ observations each (i.e., with $m$ different realizations of the
process $X$) is asymptotically as informative as observing $\mathcal{E}_0$
with $m^2n$ observations (i.e., with one realization of the process
$X$): both experiments are asymptotically equivalent to $dY_t=\sqrt
{2\sigma(t)}\,dt+m^{1/2}\delta^{1/2}n^{-1/4}\,dW_t$. Similarly, by
rescaling we can treat observations on intervals $[0,T]$ with $T>0$
fixed: observing $Y_i=X_{iT/n}+\eps_i$, $i=1,\ldots,n$, in $\mathcal{E}_0$ with $X_t=\int_0^t\sigma(s) \,dB_s$, $t\in[0,T]$, is under the
same conditions asymptotically equivalent to observing
\[
dY_u=\sqrt{2\sigma(Tu)}\, du+\delta^{1/2}T^{-1/4}n^{-1/4} \,dW_u,\qquad
u\in[0,1],
\]
or equivalently,
\[
d\tilde Y_v= \sqrt{2\sigma(v)}\, du+\delta^{1/2}(T/n)^{1/4}
\,dW_v, \qquad  v\in[0,T].
\]
Concerning the various restrictions on the smoothness $\alpha$ of the
volatility~$\sigma^2$, one might wonder whether the critical index is
$\alpha=1/2$ in view of the classical asymptotic equivalence results
[\citet{BrownLow}, \citet{Nussbaum}]. In our approach, we still face
the second-order term in \eqref{EqG3loc} and using the localized results, a much
easier globalization yields for $\alpha>1/2$ only that $\mathcal{E}_0$ is
asymptotically not less informative than observing
\[
dY_t=F(\sigma^2(t)) \,dt+ \delta^{1/2}n^{-1/4}\,dW_t,\qquad   t\in[0,1],
\]
with $F(x)=\int_1^x (y^{1/2}-2h_0^{-1})^{1/2}y^{-1}\,dy/\sqrt{8}$,
which includes a small, but nonnegligible second-order term since $h_0$
cannot tend to infinity too quickly.

On the other hand, a simple construction shows that for $\alpha<1/3$
asympto\-tic equivalence fails. In the regression model, $\mathcal{E}_0$
with $n$ observations, we cannot distinguish between $X_n(t)=\int
_0^t\sigma_n(t) \,dB_t$ with $\sigma_n^2(t)=1+n^{-\alpha}\cos(\pi
nt)$, $\norm{\sigma_n^2}_{C^\alpha}=2+n^{-\alpha}$, and standard
Brownian motion ($\sigma^2=1$) since $X_n(i/n)-X_n((i-1)/n)\sim
N(0,1/n)$ i.i.d. holds. Here, we choose the noise level $\delta
_n=n^{1/2-2\alpha}$ such that the requirement $\delta
_n^{-2}n^{-\alpha}\to0$ in Theorem \ref{ThmE0G0} holds due to
$\alpha<1/3$.

Yet, we obtain $\int_0^1(\sqrt{2\sigma_n(t)}-\sqrt{2})^2
\,dt\thicksim n^{-2\alpha}$, which shows that the signal to noise ratio
in the Gaussian shift model $\mathcal{G}_0$ with diffusion coefficient
$\delta_n^{1/2}n^{-1/4}$ is of order $n^{-2\alpha}/(\delta
_nn^{-1/2})=1$ and a Neyman--Pearson test between $\sigma_n^2$ and $1$
can distinguish both signals with a positive probability. This
different behavior for testing in $\mathcal{E}_0$ and $\mathcal{G}_0$ implies
that both models cannot be asymptotically equivalent for $\alpha<1/3$.
Note that \cit{GloterJacodI} merely require $\alpha\ge1/4$ for their
LAN-result, but our counterexample is excluded by their parametric
setting. In conclusion, the behavior in the zone $\alpha\in
[1/3,(1+\sqrt{5})/4]$ remains unexplored. If we restrict to constant
noise level $\delta$ in the regression model $\mathcal{E}_0$, then the
same argument gives a~coun\-terexample for regularity $\alpha\le1/4$.

%s8 ###
\section{Applications}\label{SecAppl}

Let us first consider the nonparametric problem of esti\-mating the spot
volatility $\sigma^2(t)$. From our asymptotic equivalence result in
Theorem~\ref{ThmE0G0} we can deduce, at least for bounded loss
functions, the usual nonparametric minimax rates, but with the number
$n$ of observations replaced by $\sqrt{n}$ provided $\sigma^2\in
C^\alpha$ for $\alpha>(1+\sqrt{5})/4$ as the mapping $\sqrt{\sigma
(t)}\mapsto\sigma^2(t)$ is a $C^\infty$-diffeomorphism for
volatilities $\sigma^2$ bounded away from zero. Since the results so
far obtained only deal with rate results, it is even simpler to use our
less informative model $\mathcal{G}_1$ or more concretely the observations
$(y_k)$ in \eqref{Eqyk} which are
independent in $\mathcal{E}_2$, centered and of variance $h^2\pi
^{-2}\sigma^2(kh)+\eps^2$. With $h=\eps$, a local (kernel or
wavelet) averaging over $\eps^{-2}\pi^2y_k^2-\pi^2$ therefore yields
rate-optimal estimators for classical pointwise or $L^p$-type loss functions.

For later use, we choose $h=\eps$ in $\mathcal{E}_2$ and propose the
simple estimator
%
%e8.1 ###
\begin{equation}\label{EqSigmaSpot}
\hat\sigma^2_b(t):=\frac{\eps}{2b}\sum_{k:\abs{k\eps-t}\le b}
(\eps^{-2}\pi^2y_k^2-\pi^2)
\end{equation}
for some bandwidth $b>0$. Since $\zeta_k^2$ is $\chi
^2(1)$-distributed, it is standard [\citet{Stone82}] to show that with
the choice $b\thicksim(\eps\log(\eps^{-1}))^{1/(2\alpha+1)}$ we
have the sup-norm risk bound
\[
\E[\norm{\hat\sigma^2_b-\sigma^2}^2_\infty]\lesssim(\eps\log
(\eps^{-1}))^{2\alpha/(2\alpha+1)},
\]
especially we shall need that $\hat\sigma^2_b$ is consistent in
sup-norm loss.

In terms of the regression experiment $\mathcal{E}_0$, we work (in an
asymptotically equivalent way) with the linear interpolation $\hat Y'$
of the observations $(Y_i)$; see the proof of Theorem \ref
{ThmRegression}. By partial integration, we can thus take for any $j,k$
%
%e8.2 ###
\begin{equation}\label{Eqyjk0}
\qquad y_{jk}^0:=-\int_0^1\Phi_{jk}(t)\hat Y''(t)\, dt=\sum_{i=1}^n
\biggl(-n\int_{(i-1)/n}^{i/n}\Phi_{jk}(t) \,dt \biggr) (Y_i-Y_{i-1}),\hspace*{-10pt}
\end{equation}
setting $Y_0:=0$. Interpreting the integral terms as weights, the
$y_{jk}^0$ are just local averages over the increments as in the
pre-averaging approach. \cit{PodolskijVetter} use Haar functions as
$\Phi_k$ (they were aware of the fact that discretized sine functions
would slightly increase the Fisher information), but they have not used
higher frequencies $j$.

Since we use the concrete coupling by linear interpolation to define
$y_{jk}^0$ in $\mathcal{E}_0$ and since convergence in total variation is
stronger than weak convergence, all asymptotics for probabilities and
weak convergence results for functionals $F((y_{jk})_{jk})$ in $\mathcal{E}_2$ remain true for $F((y_{jk}^0)_{jk})$ in $\mathcal{E}_0$, uniformly
over the parameter class. The formal argument for the latter is that
whenever $\norm{\PP_n-\QQ_n}_{\mathrm{TV}}\to0$ and $\PP_n^{X_n}\to\PP^X$
weakly for some random variables $X_n$ we have for all bounded and
continuous $g$
\[
\E_{\QQ_n}[g(X_n)]=\E_{\PP_n}[g(X_n)]+O(\norm{g}_\infty\norm{\PP
_n-\QQ_n}_{\mathrm{TV}})
\xrightarrow{n\to\infty} \E_{\PP}[g(X)].
\]
Thus, for $\alpha>1/2$, $\underline\sigma^2>0$ and $b\thicksim
(n^{-1/2}\log n)^{-1/(2\alpha+1)}$ the estimator
%
%e8.3 ###
\begin{equation}\label{EqSigmaTilde} \tilde\sigma^2_n(t):=\frac
{\delta}{2b\sqrt{n}}\sum_{k:\abs{kn^{-1/2}-t}\le b} \bigl(n\delta
^{-2}\pi^2(y_k^0)^2-\pi^2\bigr)
\end{equation}
satisfies in the regression experiment $\mathcal{E}_0$
%
%e8.4 ###
\begin{equation}\label{EqConsistentEstE0}
\lim_{n\to\infty}\inf_{\sigma^2\in\mathcal{S}(\alpha,R,\underline
{\sigma}^2)} \PP_{\sigma^2,n}\bigl(n^{\alpha/(4\alpha+2)}(\log n)^{-1}
\norm{\tilde\sigma_n^2-\sigma^2}_\infty\le R\bigr)=1.
\end{equation}

%To remain concrete, we perform a short, very rough analysis to deduce
%consistency %in sup-norm loss. For $t\in(b,1-b)$ we obtain %$\E[\hat
%$b\thicksim\eps^{1/(2\alpha+1)}$ %yields $\E[(\hat\sigma^2_b(t)-
%Since each summand in $\hat\sigma^2_b(t)$ is independent and follows a
%scaled %$\chi^2(1)$-distribution, the tails in the distribution of $
%for our choice of $b$). Consequently,

The asymptotic equivalence can be applied to construct estimators for~the
integrated volatility $\int_0^1\!\!\sigma^2(t)\,dt$ or more generally
$p$th order integrals $\int_0^1\!\!\sigma^p(t)\,dt$ using the approach
developed by \cit{IbrKhas2} for white noise models like $\mathcal{G}_0$.
In our notation, their Theorem 7.1 yields an estimator $\hat\theta
_{p,n}$ of $\int_0^1\!\!\sigma^p(t)\,dt$ in $\mathcal{G}_0$ such that
\[
\E_{\sigma^2} \biggl[ \biggl(\hat\theta_{p,n}-\int_0^1\!\!\sigma^p(t)
\,dt-\delta^{1/2}n^{-1/4}
\sqrt{2}p\int_0^1\!\!\sigma^{p-1/2}(t) \,dW_t \biggr)^2 \biggr]=o(n^{-1/2})
\]
holds uniformly over $\sigma^2\in\mathcal{S}(\alpha,R,\underline
{\sigma}^2)$ for any $\alpha,R,\underline{\sigma}^2>0$ since the
functional $\sqrt{\sigma(\bullet)}\mapsto\int_0^1\!\!\sigma^p(t)\,dt$ is
smooth on $L^2$. Their LAN-result shows that asymptotic normality with
rate $n^{-1/4}$ and variance $\delta2p^2\int_0^1\!\!\sigma^{2p-1}(t)
\,dt$ is minimax optimal. Specializing to the case $p=2$ for integrated
volatility, the asymptotic variance is $8\delta\int_0^1\!\!\sigma^3(t)
\,dt$. It should be stressed here that the existing estimation procedures
for integrated volatility are globally suboptimal for our idealized
model in the sense that their asymptotic variances involve the
integrated quarticity $\int_0^1\!\!\sigma^4(t) \,dt$ which can at most
yield optimal variance for constant values of $\sigma^2$, because
otherwise $\int_0^1\!\!\sigma^4(t) \,dt> (\int_0^1\!\!\sigma^3(t)
\,dt )^{4/3}$
follows from Jensen's inequality. The fundamental reason is that all
these estimators are based on quadratic forms of the increments
depending on
global tuning parameters, whereas optimizing weights locally permits to
attain the above efficiency bound as we shall see.

Instead of following these more abstract approaches, we use our
analysis, which is fundamentally a local likelihood approach, to
construct a simple estimator of the integrated volatility with optimal
asymptotic variance. First, we use the statistics $(y_{jk})$ in $\mathcal{E}_2$ and then transfer the results to $\mathcal{E}_0$ using $(y_{jk}^0)$
from \eqref{Eqyjk0}.

On each block $k$, we dispose in $\mathcal{E}_2$ of independent
$N(0,h^2j^{-2}\pi^{-2}\sigma^2(kh)+\eps^2)$-observations $y_{jk}$
for $j\ge1$. A maximum-likelihood estimator $\hat\sigma^2(kh)$ in
this exponential family satisfies the estimating equation
\begin{eqnarray}
\hat\sigma^2(kh)&\hspace*{3pt}=&\sum_{j\ge1} w_{jk}(\hat\sigma^2)h^{-2}j^2\pi
^2(y_{jk}^2-\eps^2),\\
\\[-28pt]
\eqntext{\mbox{where } w_{jk}(\sigma^2):=\dfrac{(\sigma^2(kh)+ h_0^{-2}\pi^2j^2)^{-2}} {\sum_{l\ge1}(\sigma^2(kh)+
h_0^{-2}\pi^2l^2)^{-2}}.}
\end{eqnarray}
This can be solved numerically, yet it is a nonconvex problem (personal
communication by J. Schmidt-Hieber). Classical MLE-theory, however,
asserts for fixed $h$, $k$ and consistent initial estimator $\tilde
\sigma_n^2(kh)$ that only one Newton step suffices to ensure
asymptotic efficiency. Because of $h\to0$ this immediate argument does
not apply here, but still gives rise to the estimator
\[
\widehat{\mathit{IV}}_\eps:=\sum_{k=0}^{h^{-1}-1} h\sum_{j\ge1}
w_{jk}(\tilde\sigma_n^2)h^{-2}j^2\pi^2(y_{jk}^2-\eps^2)
\]
of the integrated volatility $\mathit{IV}:=\int_0^1\sigma^2(t)\, dt$. Assuming
the $L^\infty$-consistency $\norm{\tilde\sigma_n^2-\sigma
^2}_\infty\to0$ in probability for the initial estimator, we assert
in $\mathcal{E}_2$ the efficiency result
\[
\eps^{-1/2}(\widehat{\mathit{IV}}_\eps-\mathit{IV})\xrightarrow{\mathcal{L}} N
\biggl(0,8\int_0^1\sigma^3(t) \,dt \biggr).
\]
To prove this, it suffices by Slutsky's lemma to show
\begin{eqnarray}\label{EqAssert1}
\qquad \eps^{-1/2}\sum_{k=0}^{h^{-1}-1} h\sum_{j\ge1} w_{jk}(\sigma
^2)h^{-2}j^2\pi^2(y_{jk}^2-\eps^2)&\xrightarrow{\mathcal{L}}& N
\biggl(0,8\int_0^1\sigma^3(t) \,dt \biggr),\\
\label{EqAssert2}
\sup_{jk}\abs{w_{jk}(\tilde\sigma_n^2)-w_{jk}(\sigma^2)}&\lesssim&
w_{jk}(\sigma^2)\norm{\tilde\sigma_n^2-\sigma^2}_\infty.
\end{eqnarray}
The second assertion \eqref{EqAssert2} follows from inserting the
Lipschitz property that $W(x):=(x+ h_0^{-2}\pi^2j^2)^{-2}$ satisfies
$\abs{W'(x)}\lesssim W(x)$, and thus $\abs{W(x)-W(y)}\lesssim
W(x)\abs{x-y}$ uniformly over $x,y\ge\underline\sigma^2>0$.

For the first assertion \eqref{EqAssert1}, note that in $\mathcal{E}_2$
the estimator $\widehat{\mathit{IV}}_\eps$ is unbiased and
\[
\Var \biggl(\sum_{j\ge1} w_{jk}(\sigma^2)h^{-2}j^2\pi
^2(y_{jk}^2-\eps^2) \biggr)=
\frac{2}{\sum_{j\ge1}(\sigma^2(kh)+h_0^{-2}\pi^2j^2)^{-2}}.
\]
We now use the identity, derived as \eqref{EqItheta},
%
%e8.5 ###
\begin{equation}\label{Eqseries}
\sum_{j\ge1}\frac{\lambda^3}{(\lambda^2+\pi^2j^2)^2}=\frac
{1+4\lambda e^{-2\lambda}-e^{-4\lambda}}{4(1-e^{-2\lambda})^2}
-\frac{1}{2\lambda}
\end{equation}
and obtain by Riemann sum approximation as $h_0\to\infty$ (with
arbitrary~speed)
\[
\eps^{-1}\Var(\widehat{\mathit{IV}}_\eps)=\sum_{k=0}^{h^{-1}-1}\frac
{2hh_0}{\sum_{j\ge1}(\sigma^2(kh)+h_0^{-2}\pi
^2j^2)^{-2}}\rightarrow8\int_0^1\sigma^3(t)\, dt.
\]
Due to the independence and Gaussianity of the $(y_{jk})$, we deduce also
\begin{eqnarray*}
&&\E \biggl[ \biggl(\sum_{j\ge1} w_{jk}(\sigma^2)h^{-2}j^2\pi
^2(y_{jk}^2-\E[y_{jk}^2]) \biggr)^4 \biggr]\\
&&\qquad \lesssim\Var \biggl(\sum_{j\ge
1} w_{jk}(\sigma^2)h^{-2}j^2\pi^2(y_{jk}^2-\eps^2) \biggr)^2
\end{eqnarray*}
such that the central limit theorem under a Lyapounov condition with
power \mbox{$p=4$} [e.g., \cit{Shiryaev}] proves assertion \eqref
{EqAssert1}, assuming $h\to0$ and $h_0\to\infty$. A~feasible
estimator is obtained by neglecting frequencies larger than some
$J=J(\eps)$:
\begin{eqnarray}\label{EqIVespJ}
\widehat{\mathit{IV}}_{\eps,J}&:=&\sum_{k=0}^{h^{-1}-1} h\sum_{j=1}^J
w_{jk}^J(\tilde\sigma_n^2)h^{-2}j^2\pi^2(y_{jk}^2-\eps^2)
\\\label{EqwjkJ}
\\[-28pt]
\eqntext{\mbox{where } w_{jk}^J(\sigma^2):=\dfrac{(\sigma^2(kh)+ h_0^{-2}\pi
^2j^2)^{-2}} {\sum_{l=1}^J(\sigma^2(kh)+ h_0^{-2}\pi
^2l^2)^{-2}}.}
\end{eqnarray}

A simple calculation yields $\E[\abs{\widehat{\mathit{IV}}_{\eps,J}-\widehat
{\mathit{IV}}_{\eps}}^2]\lesssim\eps(h_0/J)^3$ such that for $h_0/J\to0$
convergence in probability implies again by Slutsky's lemma
\[
\eps^{-1/2}(\widehat{\mathit{IV}}_{\eps,J}-\mathit{IV})\xrightarrow{\mathcal{L}} N
\biggl(0,8\int_0^1\sigma^3(t) \,dt \biggr).
\]

By the above argument, weak convergence results transfer from $\mathcal{E}_2$ to
$\mathcal{E}_0$ and we obtain the following result where we give
a concrete choice of the initial estimator, the block size $h$ and the
spectral cut-off $J$ [we just need some consistent estimator $\tilde
\sigma^2_n$, $h^{2\alpha} n^{1/2}\to0$ as well as $hn^{1/2}\to
\infty$ and $J^{-1}=o(h^{-1}n^{-1/2})$].

\begin{theorem}\label{ThmIVhat}
Let $y_{jk}^0$ for $j\ge1$, $k=0,\ldots,h^{-1}-1$ be the statistics~\eqref{Eqyjk0}
from model $\mathcal{E}_0$. For $h\thicksim n^{-1/2}\log
(n)$ and $J/\log(n)\to\infty$ consider the estima\-tor of integrated volatility
\[
\widehat{\mathit{IV}}_{n}:=\sum_{k=0}^{h^{-1}-1} h\sum_{j=1}^J
w_{jk}^J(\tilde\sigma_n^2)h^{-2}j^2\pi^2\bigl((y_{jk}^0)^2-\delta^2 n^{-1}\bigr)
\]
with weights $w_{jk}^J$ from \eqref{EqwjkJ} and the initial estimator
$\tilde\sigma_n^2$ from \eqref{EqSigmaTilde}.
Then $\widehat{\mathit{IV}}_n$ is asymptotically efficient in the sense that
\[
n^{1/4}(\widehat{\mathit{IV}}_n-\mathit{IV})\xrightarrow{\mathcal{L}} N \biggl(0,8\delta\int
_0^1\sigma^3(t) \,dt \biggr)\qquad \mbox{as }n\to\infty,
\]
provided $\sigma^2$ is strictly positive and $\alpha$-H\"{o}lder
continuous with $\alpha>1/2$.
\end{theorem}

\begin{figure}[t]

\includegraphics{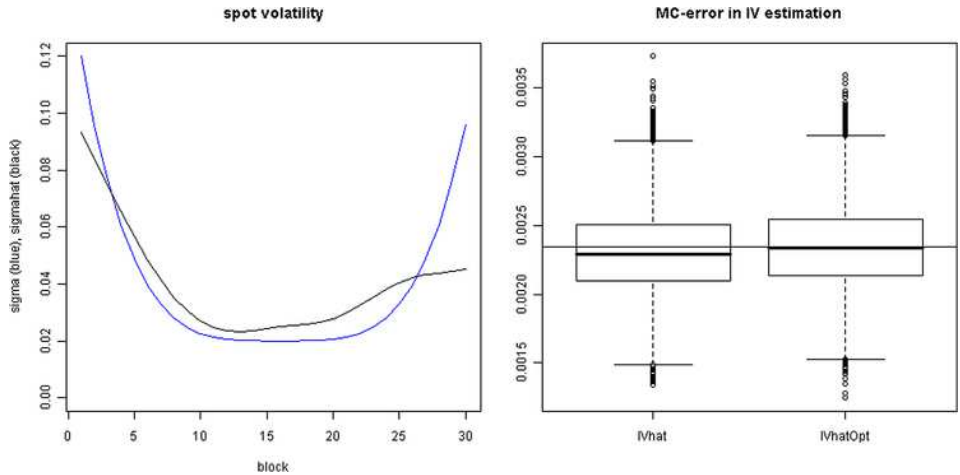}

\caption{Time-varying spot volatility and Monte Carlo error for our
estimators.}\label{Fig1}
\end{figure}

A straight-forward implementation of $\widehat{\mathit{IV}}_n$ shows a finite
sample behavior as predicted by the asymptotic results. We present
some simulation results for a situation with simplified, but realistic
model parameters. The sample size $n=30\mbox{,}000$ corresponds to roughly one
observation per second and the noise level is set to $\delta=0.01$.
The spot volatility curve $\sigma(t)=0.02+0.2(t-1/2)^4$ is
bowl-shaped, reflecting the empirical evidence of high volatility at
opening and closing. In Figure \ref{Fig1} (left) the spot volatility
and its estimate $\tilde\sigma$ on 30 blocks are presented. Instead
of \eqref{EqSigmaSpot}, we use a local-linear estimator to catch the
boundary values slightly better.
Also for the integrated volatility estimator we use $h^{-1}=30$ blocks
($h\approx6\sqrt{n}$, or expressed in real-time about 12-minute
intervals), but the estimator is quite robust to this choice.
Theoretically the maximal frequency $J$ can be as large as possible,
but due to discretization there is no more information in higher
frequencies than the block sample size. With a look at the error
analysis, we use $J:=\min(2\bar\sigma h/(\pi\delta),nh)$ with $\bar
\sigma$ denoting some upper bound on the volatility, which in our case
evaluates to $J=43$.

%f1 ###

In Figure \ref{Fig1} (right), we show the integrated volatility
estimation results obtained from 10,000 Monte Carlo iterations. The
horizontal line gives the true value $\mathit{IV}=0.0023$. The first box plot
presents the result using the weights with estimated spot volatility,
while the results with optimal oracle weights are shown in the second
box plot. We see that the estimators are practically unbiased and do
not suffer from many outliers. The empirical root mean squared error
with estimated weights is by only 5.0\% larger than the asymptotic
approximation $(8\frac{\delta}{\sqrt{n}}\int\sigma^3(t)\,dt)^{1/2}$.
With oracle weights, this reduces to $4.1\%$. An optimal procedure with
global tuning achieves asymptotically $(8\frac{\delta}{\sqrt
{n}}(\int\sigma^4(t)\,dt)^{3/4})^{1/2}$, which in our case is $19\%$
larger. Our experience with the well-established multiscale estimator
confirms this size, when oracle weights are used. Yet, it seems that
the performance of the multiscale estimator suffers significantly from
estimated weights.

Also stochastic volatility models are recovered quite well by our
implementation. The simple quadratic form of the estimator $\widehat
{\mathit{IV}}_n$ suggests that in this case a stable central limit theorem can
be derived by the usual methods. Note, however, that the analysis
cannot simply rely on our asymptotic equivalence result since $\mathcal{E}_0$ becomes non-Gaussian and, even more, Le Cam theory for stochastic
parameters (like $\sigma^2$) need to be developed. In the spirit of
\cit{Mykland}, we content ourselves with the theoretical results which
elucidate the underlying fundamental structures for the basic model and
allow straight-forward extensions to more complex models.

%s9 ###
\begin{appendix}
\section*{Appendix}\label{app}
%s9.1 ###
\subsection{Gaussian measures, Hellinger distance and Hilbert--Schmidt norm}
\label{SecPrelim}
\setcounter{equation}{0}
We gather basic facts about cylindrical Gaussian measures, the
Hellinger distance and their interplay.

Formally, we realize the white noise experiments, as $L^2$-indexed
Gaussian variables, for example, in experiment $\mathcal{E}_1$ we observe
for any $f\in L^2([0,1])$
\[
Y_f:=\scapro{f}{dY}:=\int_0^1 f(t)  \biggl(\int_0^t\sigma(s)\,dB(s)
\biggr) \,dt+\eps\int_0^1 f(t)\, dW_t.
\]
Canonically, we thus define $\PP^{\sigma,\eps}$ on the set $\Omega
=\R^{L^2([0,1])}$ with product Borel \mbox{$\sigma$-algebra} $\mathcal{F}={\mathfrak B}^{\otimes L^2([0,1])}$ (realizing a cylindrical
centered Gaussian measure). Its covariance structure is given by
\[
\E[Y_fY_g]=\scapro{Cf}{g}, \qquad f,g\in L^2([0,1]),
\]
with the covariance operator $C\dvtx L^2([0,1])\to L^2([0,1])$ given by
\[
Cf(t)=\int_0^1  \biggl(\int_0^{t\wedge u}\sigma^2(s) \,ds \biggr)f(u)
\,du+\eps^2f(t),\qquad   f\in L^2([0,1]).
\]
Note that $C$ is not trace class and thus does not define a Gaussian
measure on $L^2([0,1])$ itself.

In the construction, it suffices to prescribe $(Y_{e_m})_{m\ge1}$ for
an orthonormal basis $(e_m)_{m\ge1}$ and to set
\[
Y_f:=\sum_{m=1}^\infty\scapro{f}{e_m}Y_{e_m}.
\]
This way, we can define $\PP^{\sigma,\eps}$ equivalently on the
sequence space $\Omega=\R^{\N}$ with product $\sigma$-algebra
$\mathcal{F}={\mathfrak B}^{\otimes\N}$. This is useful when extending
results from finite dimensions.

The Hellinger distance between two probability measures $\PP$ and $\QQ
$ on $(\Omega,\mathcal{F})$ is defined as
\[
H(\PP,\QQ)= \biggl(\int_\Omega \bigl(\sqrt{p(\omega)}-\sqrt
{q(\omega)} \bigr)^2\mu(d\omega) \biggr)^{1/2},
\]
where $\mu$ denotes a dominating measure, for example, $\mu=\PP+\QQ$, and
$p$ and $q$ denote the respective densities. The total variation
distance is smaller than the Hellinger distance:
%
%e9.1 ###
\begin{equation}\label{EqHellTV}
\norm{\PP-\QQ}_{\mathrm{TV}}\le H(\PP,\QQ).
\end{equation}
The identity $H^2(\PP,\QQ)=2-2\int\sqrt{p}\sqrt{q}\,d\mu$ implies
the bound for finite or countably infinite product measures
%
%e9.2 ###
\begin{equation}\label{EqHellProd} H^2 \biggl(\bigotimes_{n}\PP
_n,\bigotimes_n\QQ_n \biggr)\le\sum_n H^2(\PP_n,\QQ_n).
\end{equation}
Moreover, the Hellinger distance is invariant under bi-measurable
bijections $T\dvtx\Omega\to\Omega'$ since with the densities $p\circ
T^{-1}$, $q\circ T^{-1}$ of the image measures~$\PP^T$ and $\QQ^T$
with respect to $\mu^T$ we have
%
%e9.3 ###
\begin{eqnarray} \label{EqHellTrafo}
H^2(\PP^T,\QQ^T)&=&\int_{\Omega'} \bigl(\sqrt{p\circ T^{-1}}-\sqrt
{q\circ T^{-1}}\bigr)^2\,d\mu^T\nonumber\\[-8pt]\\[-8pt]
&=&\int_\Omega\bigl(\sqrt{p}-\sqrt{q}\bigr)^2 \,d\mu
=H^2(\PP,\QQ).\nonumber
\end{eqnarray}
For the one-dimensional Gaussian laws $N(0,1)$ and $N(0,\sigma^2)$, we derive
\[
H^2(N(0,1),N(0,\sigma^2))=2-\sqrt{8\sigma/(\sigma^2+1)}\le2(\sigma^2-1)^2.
\]
For the multi-dimensional Gaussian laws $N(0,\Sigma_1)$ and
$N(0,\Sigma_2)$ with inverti\-ble covariance matrices $\Sigma_1,\Sigma
_2\in\R^{d\times d}$, we obtain by linear transformation and
independence, denoting by $\lambda_1,\ldots,\lambda_d$ the
eigenvalues of $\Sigma_1^{-1/2}\Sigma_2\Sigma_1^{-1/2}$:\vspace*{1pt}
\begin{eqnarray*}
H^2(N(0,\Sigma_1),N(0,\Sigma_2))&=&H^2(N(0,\Id),N(0,\Sigma
_1^{-1/2}\Sigma_2\Sigma_1^{-1/2}))\\
&\le&\sum_{k=1}^d 2(\lambda_k-1)^2.\vspace*{1pt}
\end{eqnarray*}
The last sum is nothing, but the squared Hilbert--Schmidt (or Frobenius~norm)
of $\Sigma_1^{-1/2}\Sigma_2\Sigma_1^{-1/2}-\Id$ such that\vspace*{1pt}
%
%e9.4 ###
\begin{equation}\label{EqHellHS}
H^2(N(0,\Sigma_1),N(0,\Sigma_2))\le2 \norm{\Sigma_1^{-1/2}(\Sigma
_2-\Sigma_1)\Sigma_1^{-1/2}}_{\mathrm{HS}}^2.\vspace*{1pt}
\end{equation}
Observing that \eqref{EqHellProd} and \eqref{EqHellTrafo} also apply
to Gaussian measures on the sequence space $\R^{\N}$, the bound
\eqref{EqHellHS} is also valid for (cylindrical) Gaussian measures
$N(0,\Sigma_i)$ with self-adjoint positive definite covariance
operators $\Sigma_i\dvtx L^2([0,1])\to L^2([0,1])$. %For similar
%Kullback--Leibler bounds see \cit{MunkJohannes2}.

The Hilbert--Schmidt norm of a linear operator $A\dvtx H\to H$ on any
separable real Hilbert space $H$ can be expressed by its action on an
orthonormal basis $(e_m)$ via\vspace*{-2pt}
\[
\norm{A}_{\mathrm{HS}}^2=\sum_{m,n} \scapro{Ae_m}{e_n}^2,\vspace*{1pt}
\]
which for a matrix is just the usual Frobenius norm. For self-adjoint
operators $A,B$
with $\abs{\scapro{Av}{v}}\le\abs{\scapro{Bv}{v}}$ for all $v\in H$,
we use the eigenbasis $(e_m)$ of $A$ and obtain\vspace*{1pt}
%
%e9.5 ###
\begin{equation}\label{EqHSorder}
\norm{A}_{\mathrm{HS}}^2=\sum_{m} \scapro{Ae_m}{e_m}^2\le\sum_{m,n} \scapro
{Be_m}{e_n}^2
=\norm{B}_{\mathrm{HS}}^2.\vspace*{1pt}
\end{equation}
Furthermore, it is straight-forward to see for any bounded operator $T$\vspace*{1pt}
%
%e9.6 ###
\begin{equation}\label{EqHSTrafo}
\norm{TA}_{\mathrm{HS}}\le\norm{T}\norm{A}_{\mathrm{HS}},\qquad  \norm{AT}_{\mathrm{HS}}\le
\norm{T}\norm{A}_{\mathrm{HS}}\vspace*{1pt}
\end{equation}
with the usual operator norm $\norm{T}$ of $T$. Finally, for integral
operators $Kf(x)=\int_0^1 k(x,y)f(y)\, dy$ on $L^2([0,1])$ it is well
known that\vspace*{1pt}
%
%e9.7 ###
\begin{equation}\label{EqHSL2}
\norm{K}_{\mathrm{HS}}=\norm{k}_{L^2([0,1]^2)}.\vspace*{1pt}
\end{equation}

For two Gaussian laws with different mean vectors $\mu_1,\mu_2$ and
with the same invertible covariance matrix $\Sigma$, we can similarly
use the transformation $\Sigma^{-1/2}$ and the scalar case
$H^2(N(m_1,1),N(m_2,1))=2(1-e^{-(m_1-m_2)^2/8})\le(m_1-m_2)^2/4$ to
conclude by independence\vspace*{1pt}
%
%e9.8 ###
\begin{equation}\label{EqHellmean}
H^2(N(\mu_1,\Sigma),N(\mu_2,\Sigma))\le\tfrac14 \norm{\Sigma
^{-1/2}(\mu_1-\mu_2)}^2.\vspace*{1pt}
\end{equation}
Combining \eqref{EqHellHS} and \eqref{EqHellmean}, we obtain by the
triangle inequality the bound\vspace*{-2pt}
%
%e9.9 ###
\begin{eqnarray} \label{EqHellGaussGeneral}
H^2(N(\mu_1,\Sigma_1),N(\mu_2,\Sigma_2))&\le&4\norm{\Sigma
_1^{-1/2}(\mu_1-\mu_2)}^2\nonumber\\[-9pt]\\[-9pt]
&&{}+\tfrac12\norm{\Sigma_1^{-1/2}(\Sigma_2-\Sigma_1)\Sigma_1^{-1/2}}_{\mathrm{HS}}^2.\nonumber
\end{eqnarray}

%s9.2 ###
\subsection{\texorpdfstring{Proof of Theorem \protect\ref{ThmRegression}}{Proof of Theorem 2.2.}}\label{SecProofThmRegression}

We first show that $\mathcal{E}_1$ is asymptotically at least as
informative as $\mathcal{E}_0$ for $\eps=\delta/\sqrt{n}$ and $\alpha
>0$. From $\mathcal{E}_1$ with $\eps=\delta/\sqrt{n}$, we can generate
the observations (statistics)\vspace*{-2pt}
\begin{eqnarray*}
\tilde Y_i&:=&n\int_{(2i-1)/2n}^{(2i+1)/2n}\, dY_t=n\int
_{(2i-1)/2n}^{(2i+1)/2n} X_t\,dt+\tilde\eps_i, \qquad  i=1,\ldots,n-1,
\\[-2pt]
\tilde Y_n&:=&2n\int_{(2n-1)/2n}^1 \,dY_t=2n\int_{(2n-1)/2n}^1
X_t\,dt+\tilde\eps_n,\vspace*{-2pt}
\end{eqnarray*}
with $\tilde\eps_i=n\eps(W_{(2i+1)/2n}-W_{(2i-1)/2n})\sim N(0,\delta
^2)$ and similarly $\tilde\eps_n\sim N(0,\delta^2)$, all
independent. In contrast to standard equivalence proofs, it turns out
to be essential here to take $\tilde Y_i$ as a mean symmetric around
the point $i/n$. Since $(Y_i)$ and $(\tilde Y_i)$ are defined on the
same sample space, using inequality \eqref{EqHellTV} it suffices to
prove that the Hellinger distance between the law of $(Y_i)$ and the
law of $(\tilde Y_i)$ tends to zero as $n$ tends to infinity.

For the integrated volatility function, we introduce the notation\vspace*{-2pt}
\[
a(t):=\int_0^t\sigma^2(s) \,ds,\qquad  0\le t\le1.\vspace*{-2pt}
\]
For notational convenience, we also set $a(1+s):=a(1-s)$ for $s>0$.

The covariance matrix $\Sigma^Y$ of the centered Gaussian vector
$(Y_i)$ is given~by\vspace*{-2pt}
\[
\Sigma^Y_{kl}:=\E[Y_kY_l]=a(k/n)+\delta^2\mathbf{1}(k=l),\qquad  1\le
k\le l\le n.\vspace*{-2pt}
\]
Similarly, the covariance matrix $\Sigma^{\tilde Y}$ of the centered
Gaussian vector $(\tilde Y_i)$ is given by\vspace*{-2pt}
\[
\Sigma^{\tilde Y}_{kl}:=\E[\tilde Y_k \tilde Y_l]=n\int
_{(2k-1)/2n}^{(2k+1)/2n}a(t) \,dt+\delta^2\mathbf{1}(k=l), \qquad 1\le k\le
l\le n,\vspace*{-2pt}
\]
where for $k=l=n$ we used the convention for $a(1+s)$ above. We bound
the Hellinger distance using consecutively \eqref{EqHellHS}, $\Sigma
^Y\ge\delta^2\Id$ in \eqref{EqHSorder} and \eqref{EqHellProd}, a~Taylor expansion for $a$ and treating the case $k=l=n$ by a Lipschitz
bound separately:\vspace*{-2pt}
\begin{eqnarray*}
&&H^2\bigl(\mathcal{L}(Y_i, i=1,\ldots,n),\mathcal{L}(\tilde Y_i, i=1,\ldots,n)\bigr)\\
&&\qquad \le2 \norm{(\Sigma^Y)^{-1/2}(\Sigma^Y-\Sigma^{\tilde Y})(\Sigma
^Y)^{-1/2}}_{\mathrm{HS}}^2\\
&&\qquad \le2\delta^{-4}\norm{\Sigma^{\tilde Y}-\Sigma^Y}_{\mathrm{HS}}^2\\
&&\qquad \le4\delta^{-4}\sum_{1\le k\le l\le n}  \biggl(n\int
_{(2k-1)/2n}^{(2k+1)/2n}\bigl(a(t)-a(k/n)\bigr) \,dt \biggr)^2\\
&&\qquad \le4\delta^{-4} \Biggl(O(R^2n^{-2})\\[-1pt]
&&\qquad \quad \hphantom{4\delta^{-4} \biggl(}
{}+n\sum_{k=1}^n  \biggl(n\int
_{(2k-1)/2n}^{(2k+1)/2n}\bigl(a'(k/n)(t-k/n)+O(Rn^{-1-\alpha})\bigr) \,dt
\biggr)^2 \Biggr)\\[-1pt]
&&\qquad =4\delta^{-4} \bigl(O(R^2n^{-2})+O(R^2n^{2-2-2\alpha})
\bigr)\\[-1pt]
&&\qquad =O(\delta^{-4}R^2n^{-2\alpha}).\vspace*{-2pt}
\end{eqnarray*}
Consequently, by \eqref{EqHellTV} the total-variation and thus also
the Le Cam distance between the experiments of observing $(Y_i)$ and of
observing $(\tilde Y_i)$ tends to zero for $n\to\infty$, which proves
that the white noise experiment $\mathcal{E}_1$ is asymptotically at least
as informative as the regression experiment $\mathcal{E}_0$.

To show the converse, we build from the regression experiment $\mathcal{E}_0$ a~continuous time observation by linear interpolation. To this
end, we introduce the linear $B$-splines (or hat functions)
$b_i(t)=b(t-i/n)$ with $b(t)=\min(1+nt,1-tn)\mathbf{1}_{[-1/n,1/n]}(t)$
and set\vspace*{-2pt}
\[
\hat Y'_t:=\sum_{i=1}^n Y_ib_i(t)=\sum_{i=1}^n X_{i/n}b_i(t)+\sum
_{i=1}^n\eps_i b_i(t), \qquad  t\in[0,1].\vspace*{-2pt}
\]
Note that $(\hat Y'_t)$ is a centered Gaussian process with covariance function\vspace*{-2pt}
\begin{eqnarray}
\hat c(t,s):=\E[\hat Y'_t\hat Y'_s]=\sum_{i,j=1}^n a\bigl((i\wedge
j)/n\bigr)b_i(t)b_j(s)+\delta^2\sum_{i=1}^n b_i(t)b_i(s),\nonumber \\
\eqntext{0\le t,s\le1.}\vspace*{-2pt}
\end{eqnarray}
For any $f\in L^2([0,1])$, we thus obtain\vspace*{-2pt}
\begin{eqnarray*}
\E[\scapro{f}{\hat Y'}^2]&=&\sum_{i,j=1}^n a\bigl((i\wedge j)/n\bigr)\scapro
{f}{b_i}\scapro{f}{b_j}+\delta^2\sum_{i=1}^n\scapro{f}{b_i}^2\\[-1pt]
&\le&\sum_{i,j=1}^n a\bigl((i\wedge j)/n\bigr)\scapro{f}{b_i}\scapro
{f}{b_j}+\delta^2n^{-1}\norm{f}^2,\vspace*{-2pt}
\end{eqnarray*}
because $\int nb_i=1$ yields by Jensen's inequality $\scapro
{f}{nb_i}^2\le\scapro{f^2}{nb_i}$ and we have $\sum_i b_i\le1$.
This means that the covariance operator $\hat C$ induced by the kernel
$\hat c$ is smaller than
\[
\overline C f(t):=\sum_{i,j=1}^n a\bigl((i\wedge j)/n\bigr)\scapro
{f}{b_j}b_i(t)+\delta^2n^{-1}f(t),\qquad   f\in L^2([0,1]),
\]
in the sense that $\hat C-\overline C$ is positive (semi-)definite. Now
observe that $\overline C$ is the covariance operator of the white
noise observations
%
%e9.10 ###
\begin{equation}\label{EqYbar}
d\bar Y_t= \sum_{i=1}^n X_{i/n}b_i(t) \,dt+\frac{\delta}{\sqrt
{n}}\,dW_t,\qquad   t\in[0,1].
\end{equation}
Hence, we can generate these observations from $(\hat Y'_t)$ by
randomization, that is, by adding independent, uninformative
$N(0,\overline C-\hat C)$-noise to $\hat Y'$. Now it is easy to see
that observing $\bar Y$ in \eqref{EqYbar} and $Y$ from $\mathcal{E}_1$ is
asymptotically equivalent, since in terms of the respective covariance
operators, using again \eqref{EqHellHS}, \eqref{EqHSorder} and \eqref
{EqHellProd}, the squared Hellinger distance satisfies
\begin{eqnarray*}
&&H^2(\mathcal{L}(\bar Y),\mathcal{L}(Y))\\
&&\qquad \le2\norm{(C^Y)^{-1/2}(\overline
C-C^Y)(C^Y)^{-1/2}}_{\mathrm{HS}}^2\\
&&\qquad \le2\delta^{-4}n^2\int_0^1\int_0^1  \Biggl(a(t\wedge s)-\sum
_{i,j=1}^n a\bigl((i\wedge j)/n\bigr)b_i(t)b_j(s) \Biggr)^2\,dt\,ds\\
&&\qquad = 2\delta^{-4}n^2\int_0^1\int_0^1  \Biggl(\sum_{i,j=0}^n \bigl(a(t\wedge
s)-a\bigl((i\wedge j)/n\bigr)\bigr)b_i(t)b_j(s) \Biggr)^2\,dt\,ds,
\end{eqnarray*}
where for the last line we have used $\sum_{i=0}^nb_i(t)=1$ and
$a(0)=0$. Since $b_i(t)\not=0$ can only hold when $i-\floor{nt}\in\{
0,1\}$, the $\alpha$-H\"{o}lder regularity of $\sigma^2$ implies for
$t\le s-1/n$:
\begin{eqnarray*}
&& \Biggl(\sum_{i,j=0}^n \bigl(a(t\wedge s)-a\bigl((i\wedge j)/n\bigr)\bigr)b_i(t)b_j(s)
\Biggr)^2 \\
&&\qquad =
 \Biggl(\sum_{k,l=0}^1 \bigl(a'(\floor{nt}/n)\bigl(t-(k+\floor
{nt})/n\bigr)+O(Rn^{-1-\alpha})\bigr)\\
&&\qquad \quad \hspace*{135pt}
{}\times b_{k+\floor{nt}}(t)b_{l+\floor{ns}}(s) \Biggr)^2\\
&&\qquad =O(R^2n^{-2-2\alpha})+ \Biggl(a'(\floor{nt}/n)\sum_{k=0}^1
\bigl(t-(k+\floor{nt})/n\bigr)
b_{k+\floor{nt}}(t) \Biggr)^2\\
&&\qquad =O(R^2n^{-2-2\alpha}).
\end{eqnarray*}
A symmetric argument gives the same bound for $s\le t-1/n$. For $\abs
{t-s}<1/n$, we use only the Lipschitz continuity of $a$ to obtain the
bound $O(R^2n^{-2})$. Altogether, we have found
\begin{eqnarray*}
H^2(\mathcal{L}(\bar Y),\mathcal{L}(Y))& \le&2\delta^{-4}n^2
\bigl(O(R^2n^{-2-2\alpha})+n^{-1}O(R^2n^{-2}) \bigr)
=O(\delta^{-4}R^2n^{-2\alpha}),
\end{eqnarray*}
which together with the transformation in the other direction shows
that the Le Cam distance between $\mathcal{E}_0$ and $\mathcal{E}_1$ is of
order $O(\delta^{-2}Rn^{-\alpha})$.

%s9.3 ###
\subsection{\texorpdfstring{Proof of Proposition \protect\ref{PropPiecewiseConstant}}%
{Proof of Proposition 3.2.}}
\label{SecPropPiecewiseConstant}

The main tool is Proposition \ref{PropApprox} below. Together with the
H\"{o}lder bound
\[
\abs{\sigma^2(\floor{s}_h)-\sigma^2(s)}\le R h^\alpha,\qquad   s\in[0,1],
\]
it implies that for fixed $\sigma$ the observation laws in $\mathcal{E}_1$ and $\mathcal{E}_2$ have a Hellinger distance of order $Rh^\alpha
\underline\sigma^{-3/2}\eps^{-1/2}$. By inequality \eqref
{EqHellTV}, this translates to the total variation and thus to the Le
Cam distance.

\setcounter{satz}{0}
\begin{proposition}\label{PropApprox}
For $\eps>0$ and continuous $\sigma\dvtx[0,1]\to(0,\infty)$ consider
the law $\PP^{\sigma,\eps}$ generated by
\[
dY_t= \Biggl(\int_0^t\sigma(s)\,dB(s) \Biggr) \,dt+\eps \,dW_t,\qquad   t\in[0,1],
\]
with independent Brownian motions $B$ and $W$. Then the Hellinger
distance between two laws $\PP^{\sigma_1,\eps}$ and $\PP^{\sigma
_2,\eps}$ satisfies
\[
H(\PP^{\sigma_1,\eps},\PP^{\sigma_2,\eps})\lesssim\norm{\sigma
_1^2-\sigma_2^2}_\infty \Bigl(\max_{t\in[0,1]} \sigma_1^{-3}(t)
\Bigr)\eps^{-1/2}.
\]
\end{proposition}

\begin{pf}
The covariance operator $C_\sigma$ of $\PP^{\sigma,\eps}$ is for
$f,g\in L^2([0,1])$ with antiderivatives $F,G$ satisfying $F(1)=G(1)=0$
given by
\begin{eqnarray*}
\scapro{C_\sigma f}{g} &=&\E[\scapro{f}{dY}\scapro{g}{dY}]=\E
[\scapro{f}{X}\scapro{g}{X}]+\eps^2\scapro{f}{g}\\
&=&\int FG\sigma^2+\eps^2\int fg.
\end{eqnarray*}
For covariance operators corresponding to $\sigma_1$, $\sigma_2$, we
have by twofold partial integration
\begin{eqnarray*}
\abs{\scapro{(C_{\sigma_1}-C_{\sigma_2})f}{f}}&=&\biggl|\int
_0^1\int_0^1 \int_0^{t\wedge s} (\sigma_1^2-\sigma_2^2)(u) \,du\,
f(t)f(s) \,ds \,dt\biggr|\\
&=&\biggl|\int_0^1 F(u)^2(\sigma_1^2-\sigma_2^2)(u)\, du\biggr|\\
&\le&\norm{\sigma_1^2-\sigma_2^2}_\infty\int_0^1 F(u)^2\, du
\\
&=& \norm{\sigma_1^2-\sigma_2^2}_\infty\scapro{C_{\mathrm{BM}}f}{f}
\end{eqnarray*}
with $C_{\mathrm{BM}}g(t):=\int_0^1 (t\wedge s)g(s) \,ds$, the covariance
operator of standard Brownian motion.
Using further the ordering $C_{\sigma_1}\ge\min_t\sigma
_1^2(t)C_{\mathrm{BM}}+\eps^2\Id$ and \eqref{EqHSorder}, \eqref{EqHellProd},
we obtain
\begin{eqnarray*}
&&\norm{C_{\sigma_1}^{-1/2}(C_{\sigma_2}-C_{\sigma_1})C_{\sigma
_1}^{-1/2}}_{\mathrm{HS}}\\[1pt]
&&\qquad \le\norm{\sigma_1^2-\sigma_2^2}_\infty\norm{C_{\sigma
_1}^{-1/2}C_{\mathrm{BM}}C_{\sigma_1}^{-1/2}}_{\mathrm{HS}}\\[1pt]
&&\qquad \le\norm{\sigma_1^2-\sigma_2^2}_\infty\\
&&\qquad \quad
{}\times\Bigl\|\Bigl(\min_t\sigma_1^2(t) C_{\mathrm{BM}}+\eps^2\Id\Bigr)^{-1/2}C_{\mathrm{BM}}\Bigl(\min
_t\sigma_1^2(t) C_{\mathrm{BM}}+\eps^2\Id\Bigr)^{-1/2}\Bigr\|_{\mathrm{HS}}\\
&&\qquad =\norm{\sigma_1^2-\sigma_2^2}_\infty
\norm{H(C_{\mathrm{BM}})}_{\mathrm{HS}},
\end{eqnarray*}
employing functional calculus with $H(x)=(\min_t\sigma_1^2(t)x+\eps
^2)^{-1}x$.
The spectral properties of $C_{\mathrm{BM}}$ imply that $H(C_{\mathrm{BM}})$ has
eigenfunctions $e_k(t)=\break\sqrt{2}\sin(\pi(k-1/2)t)$, $k\ge1$, with
eigenvalues $\lambda_k=\frac{4}{4\min_t\sigma_1^2(t)+(2k-1)^2\pi
^2\eps^2}$, when\-ce its Hilbert--Schmidt norm is $\norm{(\lambda
_k)}_{\ell^2}\thicksim\max_t\sigma_1^{-3/2}(t)\eps^{-1/2}$ [use
$\sum_k(s^2+k^2\eps^2)^{-2}\thicksim\eps^{-1}\times\int
(s^2+x^2)^{-2}\,dx\thicksim\eps^{-1}s^{-3}$]. This yields the result.
\end{pf}

%s9.4 ###
\subsection{\texorpdfstring{Proof of Proposition \protect\ref{PropE1G3loc}}{Proof of Proposition 5.2.}}
\label{SecProofPropE1G3loc}

We only consider the case of odd indices $k$, both cases are treated
analogously.
\cit{GramaNussbaum02} establish in their Theorem 6.1 in conjunction
with their Theorem 5.2 that $\mathcal{E}_{3,m}^{\mathrm{odd}}$ and the Gaussian
regression experiment $\mathcal{G}_{3,m}$ of observing
%
%e9.11 ###
\begin{equation}\label{EqG3m}
 \qquad\quad Y_k=v_\eps s^2(kh)+I(\sigma_0^2(kh))^{-1/2}\gamma_k,\qquad   k\in
A_m \mbox{ odd}, \gamma_k\sim N(0,1)\mbox{ i.i.d.},\hspace*{-12pt}
\end{equation}
are equivalent to experiments $\tilde\mathcal{E}_{3,m}=(\mathcal{Y},\mathcal{G},(\tilde\PP^m_{s^2})_{s^2\in C_\alpha(R)})$ and
 $\tilde\mathcal{G}_{3,m}=(\mathcal{Y},\mathcal{G},\break(\tilde\QQ^m_{s^2})_{s^2\in C_\alpha
(R)})$, respectively, on the same space $(\mathcal{Y},\mathcal{G})$ such that
%
%e9.12 ###
\begin{equation}\label{EqHellLoc}
\sup_{s^2\in C_\alpha(R)}H^2(\tilde\PP_{s^2}^m,\tilde\QQ
_{s^2}^m)\lesssim\ell^{-2\rho}
\end{equation}
holds for all $\rho<1$.

To be precise, it must be checked that the regularity conditions
(R1)--(R3) of \cit{GramaNussbaum02} are satisfied for all values $\delta
$. One complication is that in our parametric model the laws $\PP
_\theta$ and the Fisher information $I(\theta)$ depend on $h_0$ which
tends to infinity. Yet, inspecting the proofs it becomes clear that the
results remain valid if the score $\dot l=\dot l_{h_0}$ is multiplied
by $h_0^{-1/2}$ and the Fisher information accordingly by $h_0^{-1}$
and the localization is such that the parametric rate $\ell^{-1/2}$
(in our block length notation) is attained, which is ensured by our
choice in \eqref{EqParRate}. Since $I(\theta)\thicksim h_0$ is a
consequence of \eqref{EqItheta}, it remains to check conditions (R1),
(R2) of \cit{GramaNussbaum02} adjusted to our setting. Our score is
differentiable such that with $Y_j\thicksim N(0,g_j(\theta))$,
$g_j(\theta)=\theta+h_0^{-2}\pi^2j^2$
\[
\dot l_{h_0}(\theta,y)=\frac12\sum_{j\ge1}\frac{y_j^2-g_j(\theta
)}{g_j(\theta)^2},\qquad  \ddot l_{h_0}(\theta,y)=-\frac12\sum_{j\ge
1}\frac{2y_j^2-g_j(\theta)}{g_j(\theta)^3}.
\]
By the mean value theorem, (R1) requires $\E_\theta[(\ddot l(\theta
)+\frac12\dot l(\theta)^2)^2]\lesssim h_0$ (expressed in the score).
This follows here by direct moment evaluation using $\sum_{j\ge
1}g_j(\theta)^{-p}\thicksim h_0\int_0^\infty\frac{dx}{(\theta+\pi
^2x^2)^p}\thicksim h_0$ for $p>1/2$.
For (R2), we have to bound the $2\delta$-moment of $\dot l(v)\sqrt
{d\PP_v/d\PP_\theta}$ for $v$ in a neighborhood of $\theta$.
By the Cauchy--Schwarz inequality and the preceding arguments for $\dot
l$, it suffices to bound the moments of $\sqrt{d\PP_v/d\PP_\theta
}$, which are finite up to the order $\max_j\abs{1-g_j(\theta
)^2/g_j(v)^2}^{-1}$. For $v\to\theta$, this tends to infinity and
(R2) can be satisfied for any $\delta>0$. Uniform bounds are always
ensured over parameters $\theta$ bounded away from zero and infinity.

In view of the independence among the experiments
$(\mathcal{E}_{3,m}^{\mathrm{odd}})_m$ and equally among the
experiments $(\mathcal{G}_{3,m})_m$, we infer from \eqref{EqHellLoc}
and \eqref{EqHellProd}
\[
\sup_{s^2\in C_\alpha(R)}H^2\Biggl(\bigotimes_{m=1}^{(\ell h)^{-1}}\tilde\PP
_{s^2}^m, \bigotimes_{m=1}^{(\ell h)^{-1}}\tilde\QQ_{s^2}^m\Biggr)\lesssim
(\ell h)^{-1} \ell^{-2\rho} \lesssim \eps^{-1} v_\eps^2 h_0^{2\rho
}v_\eps^{4\rho}.
\]
Since we assume $h_0=o(\eps^{(1-2\alpha)/2\alpha})$, the right-hand
side tends to zero~\mbox{provided}
\[
-1+ 2\frac{\alpha}{2\alpha+1}+\frac{\rho(1-2\alpha)}{\alpha
}+\frac{4\rho\alpha}{2\alpha+1}
=\frac{\rho-\alpha}{\alpha(2\alpha+1)}>0
\]
holds. Since $\rho<1$ is arbitrary, this is always satisfied for
$\alpha<1$. In the case $\alpha=1$, we use $h_0\lesssim\eps^{-p}$
for some $p<1/2$. We have derived asymptotic equivalence between
the product experiments $\bigotimes_m\tilde\mathcal{E}_{3,m}^{\mathrm{loc}}$ and
$\bigotimes_m\tilde\mathcal{G}_{3,m}$. A fortiori, applying the \cit
{BrownLow} result, this leads to asymptotic equivalence between
observing $(y_{jk})$ in experiments $\mathcal{E}_{2,\mathrm{loc}}$ and the
corresponding Gaussian shift models of observing
%
%e9.13 ###
\begin{equation}\label{EqG2.9loc}
dY_t=I(\sigma_0^2(t))^{1/2}v_\eps s^2(t) \,dt+ (2h)^{1/2}\,dW_t,\qquad
t\in[0,1].
\end{equation}

From the explicit form \eqref{EqItheta} of the Fisher information, we
infer for $h_0\to\infty$
\[
\biggl|\frac{2\theta^{3/2}}{h_0}I(\theta)-\frac14+\frac{1}{2\theta
^{1/2} h_0}\biggr|\lesssim e^{-\underline\sigma h_0}.
\]
Consequently, by the polynomial growth of $h_0$ in $\eps^{-1}$, the
Kullback--Leibler divergence between the observation laws from
\eqref{EqG2.9loc} and the model $\mathcal{G}_{3,\mathrm{loc}}$
converges to zero. This gives the result.

%s9.5 ###
\subsection{\texorpdfstring{Proof of Proposition \protect\ref{PropE1locyjk}}{Proof of Proposition 6.1.}}
\label{ProofPropE1locyjk}

Since the observations $y_{jk}$ for $j\ge1$ are the same in $\mathcal{Y}$
and $\tilde\mathcal{Y}$, we can work conditionally on those. Moreover, it
suffices to consider only the event $\Omega_\eps:=\{\norm{\hat
\sigma_\eps^2-\sigma^2}_\infty\le Rv_\eps\}$ because the squared
Hellinger distance satisfies by conditioning and restriction to $\Omega
_\eps$ (with density functions $f$ and further obvious notation)
\begin{eqnarray*}
H^2(\mathcal{L}(\mathcal{Y}),\mathcal{L}(\tilde\mathcal{Y}))
&=&\int \bigl(\sqrt{f_{\mathcal{Y}|(y_{jk})_{j\ge1,k}}f_{(y_{jk})_{j\ge
1,k}}}-\sqrt{f_{\tilde\mathcal{Y}|(y_{jk})_{j\ge1,k}}f_{(y_{jk})_{j\ge
1,k}}} \bigr)^2\\
&=&\E\bigl[H^2\bigl(\mathcal{L}\bigl((y_{0k})_k | (y_{jk})_{j\ge1,k}\bigr),\mathcal{L}\bigl((\tilde
y_{0k})_k | (y_{jk})_{j\ge1,k}\bigr)\bigr)\bigr]\\
&\le&\E\bigl[H^2(\mathcal{L}\bigl((y_{0k})_k | (y_{jk})_{j\ge1,k}\bigr),\mathcal{L}\bigl((\tilde y_{0k})_k | (y_{jk})_{j\ge1,k}\bigr))\mathbf{1}_{\Omega_\eps}\bigr]\\
&&{} +2\PP(\Omega_\eps^\complement)
\end{eqnarray*}
with $\PP(\Omega_\eps^\complement)\to0$. Conditional on
$(y_{jk})_{j\ge1,k}$, both laws are Gaussian, $(y_{0,k})_k$ has
mean $\mu$ with
\begin{eqnarray*}
\mu_0&=&2\sum_{j\ge1}\frac{\Var(\beta_{jk})}{\Var
(y_{jk})}y_{j0},\\
\mu_k&=&\sum_{j\ge1} \biggl(\frac{\Var(\beta
_{j,k-1})}{\Var(y_{j,k-1})}(-1)^{j+1}y_{j,k-1}
+\frac{\Var(\beta_{j,k-1})}{\Var(y_{j,k-1})}y_{jk} \biggr)
\end{eqnarray*}
for $k\ge1$ and covariance matrix $\Sigma$ with
\[
\Sigma_{k,k'}=
\cases{
\displaystyle c_k\eps^2\sum_{j\ge1} \biggl(\frac{\Var(\beta_{j,k-1})}{\Var
(y_{j,k-1})}+\frac{\Var(\beta_{jk})}{\Var(y_{jk})} \biggr) +\eps^2,&
\quad if $ k'=k$,\cr
\displaystyle c_{k\wedge k'}\eps^2\sum_{j\ge1}(-1)^{j+1} \frac{\eps^2\Var(\beta
_{j,k-1})}{\Var(y_{j,k-1})}
-\frac{\eps^2}{2},&\quad  if $k'=k\pm1$,\cr
0,&\quad otherwise,
}
%
%  k,k'=0,\ldots,h^{-1}-1,
\]
where $c_k:=1\vee(2-k)\in\{1,2\}$. Conditional mean $\tilde\mu$ and
covariance matrix $\tilde\Sigma$ of $(\tilde y_{0k})_k$ have the same
representation, but replacing $\Var$ each time by $\Var_\eps$,
compare~\eqref{EqVareps}.

%From the tri-diagonal structure of $\Sigma$ and from
From $\frac{\Var(\beta_{jk})}{\Var(y_{jk})}=(1+h_0^{-2}\pi
^2j^2\sigma^2(kh))^{-1}$, we infer for $h_0\to\infty$ by Riemann sum
approximation
\begin{eqnarray*}
\sum_{j\ge1} \biggl(\frac{\Var(\beta_{j,k-1})}{\Var
(y_{j,k-1})}+\frac{\Var(\beta_{jk})}{\Var(y_{jk})} \biggr) &\thicksim&
\sum_{j\ge1}\frac{1}{1+j^2h_0^{-2}}\thicksim h_0,\qquad   h_0\to\infty
,\\[-2pt]
\Biggl|\sum_{j\ge1}(-1)^{j+1}
\frac{\Var(\beta_{j,k-1})}{\Var(y_{j,k-1})}\Biggr|
&\thicksim&\sum_{j\ge1}\frac
{2jh_0^{-2}}{(1+(2j)^2h_0^{-2})(1+(2j+1)^2h_0^{-2})}\thicksim1.
\end{eqnarray*}
Hence, $\Sigma$ is a matrix with entries of order $\eps^2h_0$ on the
main diagonal and entries of order $\eps^2$ on the two adjacent
diagonals. A simple Cauchy--Schwarz argument therefore shows $\scapro
{\Sigma v}{v}\gtrsim(\eps^2h_0-\eps^2)\norm{v}^2\thicksim\eps
^2h_0\norm{v}^2$ for $h_0\to\infty$
which implies $\Sigma\gtrsim\eps h\Id$ in matrix order.
Combining this with the Hellinger bound \eqref{EqHellGaussGeneral}, we
arrive at the estimate
\begin{eqnarray*}
&&\E\bigl[H^2\bigl(\mathcal{L}\bigl((y_{0k})_k | (y_{jk})_{j\ge1,k}\bigr),\mathcal{L}\bigl((\tilde
y_{0k})_k | (y_{jk})_{j\ge1,k}\bigr)\bigr)\bigr]\\
&&\qquad \lesssim\E \biggl[\frac{\norm{\mu-\tilde\mu}^2}{\eps h}
\biggr]+\frac{\norm{\Sigma-\tilde\Sigma}_{\mathrm{HS}}^2} {\eps^2h^2}\\
&&\qquad \lesssim
\sum_{j\ge1,k}  \biggl(\frac{\Var(\beta_{jk})}{\Var(y_{jk})}
-\frac{\Var_\eps(\beta_{jk})}{\Var_\eps(y_{jk})} \biggr)^2\frac
{\Var(y_{jk})}{\eps h}\\[-2pt]
&&\qquad \quad {} +
\sum_{j\ge1,k}  \biggl(\frac{\eps^2\Var(\beta_{jk})}{\Var(y_{jk})}
-\frac{\eps^2\Var_\eps(\beta_{jk})}{\Var_\eps(y_{jk})}
\biggr)^2\eps^{-2}h^{-2}.
%&\lesssim
%-\frac{\eps^2\Var_\eps(\beta_{jk})}{\Var_\eps(y_{jk})} )^2
\end{eqnarray*}
The function $G(z):=\frac{\norm{\Phi_{jk}}^2z}{\norm{\Phi
_{jk}}^2z+\eps^2}$ has derivative $G'(z)=\frac{\norm{\Phi
_{jk}}^2\eps^2}{(\norm{\Phi_{jk}}^2z+\eps^2)^2}$ and thus satisfies
uniformly over all $z$ bounded away from zero
$\abs{G(w)-G(z)}\lesssim\frac{\norm{\Phi_{jk}}^2\eps^2\abs
{w-z}}{(\norm{\Phi_{jk}}^2+\eps^2)^2}$.\vspace{2pt} Inserting $\abs{\sigma
^2-\sigma_0^2}\lesssim v_\eps$ and $\norm{\Phi_{jk}}\thicksim h/j$,
we thus find the uniform bound on $\Omega_\eps$
\[
\biggl(\frac{\Var(\beta_{jk})}{\Var(y_{jk})}
-\frac{\Var_\eps(\beta_{jk})}{\Var_\eps(y_{jk})} \biggr)^2\lesssim
\frac{v_\eps^2\eps^4h^4/j^4}{(\eps^2+h^2/j^2)^4}\thicksim v_\eps
^2\min(h_0/j,j/h_0)^4.
\]
Putting the estimates together, we arrive at\vspace{-2pt}
\begin{eqnarray*}
H^2(\mathcal{L}(\mathcal{Y}),\mathcal{L}(\tilde\mathcal{Y}))&\lesssim&
v_\eps^2\sum_{j\ge1,k} \min(h_0/j,j/h_0)^4 \biggl(\frac
{1+h_0^2/j^2}{h_0}+\frac{1}
{h_0^2} \biggr)+\PP(\Omega_\eps^\complement)\\
&\le&2v_\eps^2h^{-1}\sum_{j\ge1}\min(h_0/j,j/h_0)^2h_0^{-1}+\PP
(\Omega_\eps^\complement)\\
&\thicksim& v_\eps^2h_0^{-1}\eps^{-1}+\PP(\Omega_\eps^\complement)
\end{eqnarray*}
such that the Hellinger distance tends to zero uniformly if
$h_0^{-1}v_\eps^2=o(\eps)$, which is ensured by our choice of $h_0$.
This implies asymptotic equivalence of observing $\mathcal Y$ and $\tilde
\mathcal{Y}$ and thus of experiment $\mathcal{E}_2$ and of just observing
$(y_{jk})_{j\ge1,k}$ in $\mathcal{E}_2$. By independence, the latter is
equivalent to $\mathcal{E}_{2,\mathrm{odd}}\otimes\mathcal{E}_{2,\mathrm{even}}$.
\end{appendix}

\section*{Acknowledgments}
I am grateful to Marc Hoffmann, Mark Podolskij and Johannes
Schmidt-Hieber for very
useful discussions and to three referees and an associate editor for
their very careful reading and helpful comments.

%suskaldyti doi

\printaddresses

\end{document}